\documentclass[12 pt]{article}
\usepackage [english]{babel}
\usepackage {pgfpict2e}
\usepackage{amssymb}
\usepackage{amsthm}
\usepackage{amsmath}
\usepackage{setspace}
\usepackage{euscript}
\usepackage{mathrsfs}
\usepackage{graphicx}
\usepackage{textcomp}
\usepackage{hyperref}
\usepackage{tikz}
\usetikzlibrary{calc}

\newtheorem{Claim}{Claim}
\newtheorem{Conjecture}{Conjecture}
\newtheorem{Def}{Definition}
\newtheorem{conseq}{Corollary}
\newtheorem{Lemma}{Lemma}
\newtheorem{Remark}{Remark}

\newtheorem{Theorem}{Theorem}

\newcommand{\FF}[2]{#1^{\underline{#2}}}
\newcommand{\NPK}[2]{\genfrac{\{}{\}}{0pt}{}{#1}{#2}}
\newcommand{\sltwo}{{\mathfrak{sl}_2}}
\newcommand{\ti}{\widetilde}

\title{Values of the $\mathfrak{sl}_2$ weight system on complete bipartite graphs.}
\author{P.~Filippova \thanks{National Research University Higher School of Economics,
e-mail: apoly38@gmail.com. 
The publication was prepared within the framework of the Academic Fund Program at HSE University in 2020–2021 (grant \textnumero 20-04-010) and within the framework of the Russian Academic Excellence Project  ``5-100''.
}}
\begin{document}

\maketitle

\begin{abstract}
A weight system is a function on chord diagrams that satisfies the so-called four-term relations. Vassiliev's theory of finite-order knot invariants describes these invariants in terms of weight systems. In particular, there is a weight system corresponding to the colored Jones polynomial. This weight system can be easily defined in terms of the Lie algebra $\mathfrak{sl}_2$, but this definition is too cumbersome from the computational point of view, so that the values of this weight system are known only for some limited classes of chord diagrams. 

In the present paper we give a formula for the values of the $\mathfrak{sl}_2$ weight system for a class of chord diagrams whose intersection graphs are complete bipartite graphs with no more than three vertices in one of the parts. 

Our main computational tool is the Chmutov--Varchenko reccurence relation. Furthermore, complete bipartite graphs with no more than three vertices in one of the parts generate Hopf subalgebras of the Hopf algebra of graphs, and we deduce formulas for the projection onto the subspace of primitive elements along the subspace of decomposable elements in these subalgebras. We compute the values of the $\mathfrak{sl}_2$ weight system for the projections of chord diagrams with such intersection graphs. Our results confirm certain conjectures due to S.K.Lando on the values of the weight system $\mathfrak{sl}_2$ at the projections of chord diagrams on the space of primitive elements.
\end{abstract}
\tableofcontents
\section{Introduction} 

Finite-order knot invariants, which were introduced in~\cite{V90} by Vassiliev about 1990, can be expressed in terms of weight  systems, that is, functions on chord diagrams satisfying the so-called Vassiliev four-term relations. In the paper~\cite{Ko}, Kontsevitch proved that over a field of characteristic zero every weight system corresponds to some finite-order invariant.

There are multiple approaches to constructing weight systems. In particular, Bar-Natan~\cite{BN} and Kontsevitch~\cite{Ko} suggested a construction of a weight system from a finite-dimensional Lie algebra endowed with an invariant nondegenerate bilinear form. The simplest case of this construction is the $\sltwo$ weight system, which is constructed from the Lie algebra $\sltwo$. Its values lie in the center of the universal enveloping algebra of $\sltwo$. The center is isomorphic to the ring of polynomials in one variable (the Casimir element). 
In Vassiliev's approach this weight system corresponds to a well-known knot ivariant, the colored Jones polynomial. 
The value of the $\sltwo$ weight system at a chord diagram of order~$n$ is a monic polynomial of degree~$n$ in the Casimir element.

The $\sltwo$ weight system was studied in many papers. Despite the simplicity of the definition of this weight system, it is difficult to compute its value at a chord diagram by using this definition, because doing this involves computations in a noncommutative algebra.
The Chmutov--Varchenko recurrence relations~\cite{bib:ChV} significantly simplify these computations; however, using them in computations is laborious as well, and the explicit values of the $\sltwo$ weight system are known only for chord diagrams of low order and for a small number of simple series of chord diagrams. 
In particular, the values of the $\sltwo$ weight system at chord diagrams with complete intersection graph are unknown. The conjecture of Lando about the form of the corresponding polynomials has been proved only for the linear terms of these polynomials~\cite{Bigeni}.

The Chmutov--Lando theorem~\cite{ChL} states that the value of the $\sltwo$ weight system at a chord diagram depends only on the intersection graph of this chord diagram; i.e., if two chord diagrams have isomorphic intersecton graphs, then the values of the weight system at these chord diagrams coincide. 
This raises the following natural question, which was asked by Lando in~\cite{ChL}: Is it possible to extend this weight system to a polynomial graph invariant satisfying the four-term relations for graphs? E.~Krasil'nikov showed that such an extension exists and is unique for all graphs with $n\le8$ vertices, but nothing is known about such an extension for graphs with more than eight vertices. One of the possible approaches to answer this question is to \emph{define} a polynomial invariant for arbitrary graphs and then show that this invariant satisfies the four-term relations for graphs and coincides with the $\sltwo$ weight system on the intersection graphs. 
To make this possible, it is necessary to have enough examples of the explicit values of the $\sltwo$ weight system for various graph families. In our present paper we compute such values at the complete bipartite graphs such that the number of vertices in one of the parts is at most three.

The quotient space of the vector space of chord diagrams by the four-term relations carries the structure of a connected graded commutative cocommutative Hopf algebra~\cite{Ko}.
The same is true for the vector space of graphs.
Further, the complete bipartite graphs generate a Hopf subalgebra in the Hopf algebra of graphs, and for any $l = 0,1,2,\ldots ,$ this Hopf subalgebra contains the Hopf subalgebra generated by complete bipartite graphs such that one of the parts contains no more than $l$ vertices.

According to the Milnor--Moore theorem~\cite{MM},
each of these Hopf algebras is generated by its primitive elements
and is a polynomial Hopf algebra in its primitive elements.
As a consequence, in each of these Hopf algebras there is a well-defined projection onto
the subspace of primitive elements along the space of decomposable elements.
There is a universal formula which expresses this projection as a logarithm 
of the identity homomorphism~(\cite{L97,S94}).
However, this formula is cumbersome, and we construct
its compact forms for the Hopf algebras of complete bipartite graphs with sizes of one of the parts 
at most $1,2,$~and~$3$. Then we use these forms to explicitly 
compute  of the values of the $\sltwo$ weight system at
the projections on the subspace of primitive elements of the complete bipartite graphs with size of one of the parts
at most $l = 1,2,3$. It turns out that these values are polynomials of degree at most $l$. Thus, our computations confirm another 
conjecture of Lando, which states that 
if $C$ is a chord diagram such that the circumference (that is, the length of the longest cycle) of its intersection graph is at most $2l, l\ge1$,
then  the value of the  $\sltwo$ weight system at its projection on the subspace of primitive elements
is a polynomial of degree at most~$l$.

Note that the value of a multiplicative graph invariant can be uniquely reconstructed from
the values of this invariant at the projections of graphs on the subspace of primitive elements.
What is more, this invariant is often significantly simplified under the projection.
For instance, it is known
 that the value of the $\sltwo$ weight system at the projection of a chord diagram of order~$n$ on the subspace of primitive elements is a polynomial of degree less than or equal to $[\frac{n}2]$, 
i.e., the degree of the polynomial is at least halved.
In~\cite{KLMR} the coefficient of the term of degree~$\frac{n}2$ was extended to  arbitrary graphs satisfying the four-term relations for graphs. This gives the hope for the existence of such an extension for the other coefficients of the polynomial. 
(In~\cite{LZ17}, an extension of this coefficient for the more general case of binary delta-matroids was given.)

This paper is organized as follows.
In Section~\ref{sec:HA} we give definitions of the Hopf algebras of chord diagrams and graphs and of their subalgebras that we are interested in. Our main new result here is explicit formulas for the projections of complete bipartite graphs with size of one of the parts at most three on the subspace of primitive elements.
These formulas may be useful for computing  such projections for other graph invariants.

In Section~\ref{sec:sl2} we compute explicitly the values of the $\sltwo$ weight system
at the chord diagrams whose intersection graph is a complete bipartite graph with size of one 
of the parts at most three.
Using this result and the projection formula from Section~\ref{sec:HA}, we compute the values of the $\sltwo$ weight system at the projections of such complete bipartite graphs on the subspace of primitive elements.

We follow the approach of~\cite{ZvL}; see also~\cite{cd}.

\section{Hopf algebras of graphs and chord diagrams}\label{sec:HA}

In this section we define the Hopf algebra of graphs and the Hopf algebra of chord diagrams modulo the four-term relations. We also define Hopf subalgebras generated by complete bipartite graphs in the Hopf algebra of graphs. Using the universal formula for the projection onto the subspace of primitive elements, we derive formulas for the projection onto the subspace of primitive  elements in the Hopf algebras of complete bipartite graphs with no more than $l$, $l = 1,2,3$ vertices in one of the parts. 

A \emph{counital coassociative coalgebra} over a field $\mathbb K$ is a vector space $C$ over $\mathbb K$ together with $\mathbb K$-linear maps 
\begin{align*}
\mu& \colon C \to C \otimes C\\
\varepsilon& \colon C \to \mathbb K
\end{align*}
such that
\begin{align*}
(id_C \otimes \mu)\circ \mu = (\mu \otimes id_C)\circ \mu, \\
(id_C \otimes \varepsilon)\circ \mu = id_C = (\varepsilon \otimes id_C)\circ \mu.
\end{align*}

A \emph{bialgebra} over a field $\mathbb K$ is a vector space $B$ over $\mathbb K$ endowed with two structures, that of a unital associative algebra over $\mathbb K$ (with multiplication $m$ and unit $\eta$) and that of a counital coassociative coalgebra (with comultiplication $\mu$ and counit $\varepsilon$) such that
\begin{align*}
\mu \circ m = (m \otimes \mu) \circ (id \otimes \tau \otimes id) \circ (\mu \otimes \mu), \\
\varepsilon \otimes \varepsilon  = \varepsilon \circ m, \\
\eta \otimes \eta = \mu \circ \eta, \\
id = \varepsilon \circ \nu.
\end{align*}

Here by $\tau \colon B \otimes B \to B$ we denote the linear map defined by $\tau (x \otimes y) = y \otimes x, x, y \in B.$

A \emph{Hopf algebra} over a field $\mathbb K$ is an associative, unital, coassociative, and counital bialgebra $H$ together with an \emph{antipode}, that is, a $\mathbb K$-linear map  $S \colon H \to H$ such that (in the above notation)
\begin{align*}
m \circ (S \otimes id) \circ \mu = \eta \circ \varepsilon = m \circ (id \otimes S) \circ \mu.
\end{align*}

Further on we assume that the characteristic of the ground field $\mathbb K$ is zero.

\subsection{Hopf algebras of graphs and chord diagrams}

In this paper by a graph we mean an isomorphism class of finite simple graphs (i.e. finite graphs with no loops and multiple edges). Formal linear combinations of graphs form a vector space, which is graded by the number of graph vertices

We define the \emph{product of two graphs} $G_1$ and $G_2$ as their disjoint union: $G_1 G_2 := G_1 \sqcup G_2$.

This multiplication is extended to the vector space of graphs by linearity. It preserves the grading. Thus, this vector space is endowed with the structure of a graded algebra.

By $V(G)$ we denote the vertex set of a graph $G$.

Any subset $U\subset V(G)$ induces a subraph of $G$. We denote such a subgraph by $G|_U$. 
We define the \emph{comultiplication} $\mu$ which acts on a graph~$G$  as
\begin{equation*}
\mu(G):= \sum_{U\subset V(G)} G|_{U}\otimes G|_{V(G)\setminus U}.\label{graphs_coprod}
\end{equation*}

Both multiplication and comultiplication are extended to the vector space of graphs by linearity and preserve the grading. Thus, this vector space is endowed with both the structure of a graded algebra and that of a graded coalgebra. Moreover, the following assertion holds.

\begin{Claim}
The multiplication and comultiplication defined above, together with the naturally defined unit, counit, and antipode, turn the vector space of graphs into a Hopf algebra.
\end{Claim}

This construction was introduced in~\cite{JR}.

By $\EuScript G$ we denote the Hopf algebra of graphs. The set of all graphs with $n$ vertices generate a vector subspace $\EuScript G_n$ in $\EuScript G$. Thus, 

\begin{equation*}
\EuScript G = \EuScript G_0\oplus \EuScript G_1 \oplus \EuScript G_2 \oplus \ldots
\end{equation*}

Let $A$ and $B$ be two vertices of a graph $G$. By $G'_{AB}$ we denote the graph obtained from $G$ by changing the adjacency between the vertices $A$ and $B$ in $\Gamma$, i.e., by deleting the edge $AB$ if it exists and adding the edge $AB$ otherwise. By $\tilde G_{AB}$ we denote the graph obtained from $G$ by changing the adjacency with $A$ of each vertex in $V(G)\setminus \{A,B \}$ joined with $B$. 
A \emph{four-term element in the space of graphs} is a linear combination
\begin{equation*}
G - G'_{AB} - \tilde G_{AB} + \tilde G_{AB}',
\end{equation*}
Note that all graphs in a four-term element have the same number of vertices. 

Let $\EuScript F_n$ denote the quotient space of $\EuScript G_n$ by the subspace spanned by the four-term elements containing $n$-vertex graphs. The graded Hopf algebra structure on $\EuScript G$ induces a graded Hopf algebra structure on the space $\EuScript F$:
\begin{equation*}
\EuScript F = \EuScript F_0\oplus \EuScript F_1 \oplus \EuScript F_2 \oplus \ldots
\end{equation*}

Our description of the Hopf algebra structure on graphs follows \cite{Lando}.


\begin{Def}
A \emph{chord diagram} of order $n$ (a chord diagram with $n$ chords) is an oriented circle together with $2n$ pairwise distinct points splitted into $n$ disjoint pairs considered up to orientation-preserving diffeomorphisms of the circle. 
\end{Def}

We connect the points belonging to the same pair by a segment of a line or of a curve, called a \emph{chord}. (The shape of a chord is irrelevant, but it must have no common points with the circle except its endpoints.) 

The vector space spanned by the chord diagrams is graded. Each componet is spanned by diagrams of the same order. 

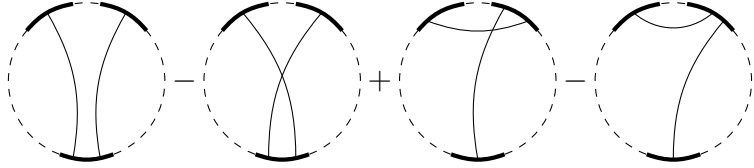
\begin{figure}
\begin{center}
\begin{tikzpicture}[scale = 1.3]
\draw[dashed]  (1,3.5) circle [radius=0.8];
\node  at (2,3.5) {$-$};
\draw[dashed]  (3,3.5) circle [radius=0.8];
\node  at (4,3.5) {$+$};
\draw[dashed]  (5,3.5) circle [radius=0.8];
\node  at (6,3.5) {$-$};
\draw[dashed]  (7,3.5) circle [radius=0.8];

\foreach \n in {0,2,...,6}{
\node(center) at (1 + \n,3.5) {};
\path ($(center) + (40:0.8)$) edge[bend right =20, ultra thick] ($(center) + (80:0.8)$);
\path ($(center) + (100:0.8)$) edge[bend right =20, ultra thick] ($(center) + (140:0.8)$);
\path ($(center) + (250:0.8)$) edge[bend right =20, ultra thick] ($(center) + (290:0.8)$);
}
\path  ($(1,3.5) + (120:0.8)$) edge[bend left =20] ($(1,3.5) + (260:0.8)$);
\path  ($(1,3.5) + (60:0.8)$) edge[bend right =20] ($(1,3.5) + (280:0.8)$);
\path  ($(3,3.5) + (120:0.8)$) edge[bend left =20] ($(3,3.5) + (280:0.8)$);
\path  ($(3,3.5) + (60:0.8)$) edge[bend right =20] ($(3,3.5) + (260:0.8)$);
\path  ($(5,3.5) + (130:0.8)$) edge[bend right =20] ($(5,3.5) + (50:0.8)$);
\path  ($(5,3.5) + (70:0.8)$) edge[bend right =20] ($(5,3.5) + (270:0.8)$);
\path  ($(7,3.5) + (120:0.8)$) edge[bend right =40] ($(7,3.5) + (60:0.8)$);
\path  ($(7,3.5) + (50:0.8)$) edge[bend right =20] ($(7,3.5) + (270:0.8)$);
\end{tikzpicture}
\caption{4-term element in the space of chord diagrams} \label{pict:4term}
\end{center}
\end{figure}

A \emph{four-term element} in the space of chord diagrams  is the linear combination of diagrams which is shown in Fig.~\ref{pict:4term}.
The four diagrams in this combinations contain the same set of chords in addition to those shown in the figure, but the endpoints of all these chords must belong to the dashed arcs.
 
Equating four-term elements in spaces of graphs and chord diagrams to zero, we obtain \emph{four-term relations} in the corresponding spaces.

\begin{Def}
An \emph{arc diagram} of order $n$ is an oriented line together with $2n$ pairwise distinct points splitted into $n$ disjoint pairs,  which is considered up to orientation-preserving diffeomorphisms of the line. 
\end{Def}

Each of these $n$ pairs of points is shown as an arc joining these points and lying in the upper half-plane. 

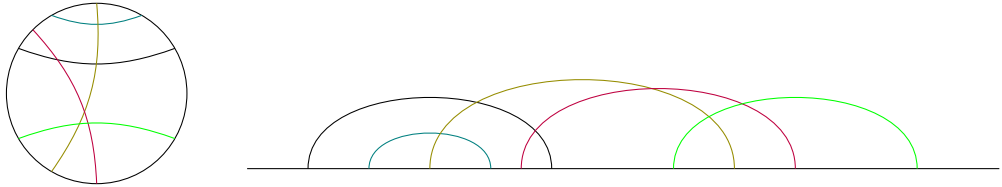
\begin{figure}
\begin{center}
\begin{tikzpicture}
\draw  (2,3.5) circle [radius=1.2];
\draw (4,2.5)--(14,2.5);
\node(center) at (2,3.5) {};
\node(origin) at (4,2.5) {};
\path  ($(center) + (30:1.2)$) edge[bend left =20, color = black] ($(center) + (150:1.2)$);
\path  ($(origin) + (30*0.027,0)$) edge[bend left =90, color = black] ($(origin) + (150*0.027,0)$);
\node at ($(origin) + (30*0.027,-0.25)$) {};
\node at ($(origin) + (160*0.027,-0.25)$) {};
\path  ($(center) + (60:1.2)$) edge[bend left =20, color = teal] ($(center) + (120:1.2)$);
\path  ($(origin) + (60*0.027,0)$) edge[bend left =90, color = teal] ($(origin) + (120*0.027,0)$);
\node at ($(origin) + (60*0.027,-0.25)$) {};
\node at ($(origin) + (118*0.027,-0.25)$) {};
\path  ($(center) + (90:1.2)$) edge[bend left =20, color = olive] ($(center) + (240:1.2)$);
\path  ($(origin) + (90*0.027,0)$) edge[bend left =90, color = olive] ($(origin) + (240*0.027,0)$);
\node at ($(origin) + (90*0.027,-0.25)$) {};
\node at ($(origin) + (240*0.027,-0.25)$) {};
\path  ($(center) + (135:1.2)$) edge[bend left =20, color = purple] ($(center) + (270:1.2)$);
\path  ($(origin) + (135*0.027,0)$) edge[bend left =90, color = purple] ($(origin) + (270*0.027,0)$);
\node at ($(origin) + (141*0.027,-0.25)$) {};
\node at ($(origin) + (270*0.027,-0.25)$) {};
\path  ($(center) + (210:1.2)$) edge[bend left =20, color = green] ($(center) + (330:1.2)$);
\path  ($(origin) + (210*0.027,0)$) edge[bend left =90, color = green] ($(origin) + (330*0.027,0)$);
\node at ($(origin) + (210*0.027,-0.25)$) {};
\node at ($(origin) + (330*0.027,-0.25)$) {};
\end{tikzpicture}
\caption{An example of an arc representation of a chord diagram.\label{pic:chord-arc}}
\end{center}
\end{figure}

Choosing a point on a chord diagram (different from the endpoints of all chords) and cutting the chord diagram at this point, we obtain an arc diagram, or an \emph{arc representation,} of this chord diagram (see Fig.~\ref{pic:chord-arc}). 
A chord diagram may have up to $2n$ different arc representations, while an arc diagram uniquely determines the corresponding chord diagram.

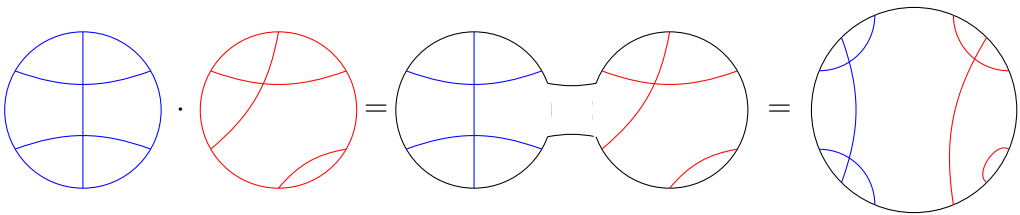
\begin{figure}
\begin{center}
\begin{tikzpicture}[scale = 1.3]
\draw[color = blue]  (1,3.5) circle [radius=0.8];
\node  at (2,3.5) {$\cdot$};
\draw[color = red]  (3,3.5) circle [radius=0.8];
\node  at (4,3.5) {$=$};
\draw  (5,3.5) circle [radius=0.8];
\draw  (7,3.5) circle [radius=0.8];
\node  at (8.125,3.5) {$=$};
\draw  (9.5,3.5) circle [radius=1.05];
\path  ($(5,3.5) + (20:0.8)$) edge[bend left =23.5, color = white, ultra thick] ($(5,3.5) + (340:0.8)$);
\path  ($(5,3.5) + (340:0.8)$) edge[bend left = 10, color = black] ($(7,3.5) + (200:0.8)$);
\path  ($(7,3.5) + (200:0.8)$) edge[bend left =23.5, color = white, ultra thick] ($(7,3.5) + (160:0.8)$);
\path   ($(7,3.5) + (160:0.8)$) edge[bend left = 10, color = black]   ($(5,3.5) + (20:0.8)$);
\path  ($(1,3.5) + (150:0.8)$) edge[bend right =20, color = blue] ($(1,3.5) + (30:0.8)$);
\path  ($(1,3.5) + (210:0.8)$) edge[bend left =20, color = blue] ($(1,3.5) + (330:0.8)$);
\path  ($(1,3.5) + (90:0.8)$) edge[color = blue] ($(1,3.5) + (270:0.8)$);
\path  ($(3,3.5) + (-150:0.8)$) edge[bend right = 20, color = red] ($(3,3.5) + (-270:0.8)$);
\path  ($(3,3.5) + (-210:0.8)$) edge[bend right = 20, color = red] ($(3,3.5) + (-330:0.8)$);
\path  ($(3,3.5) + (-90:0.8)$) edge[bend left = 20, color = red] ($(3,3.5) + (-30:0.8)$);
\path  ($(5,3.5) + (150:0.8)$) edge[bend right =20, color = blue] ($(5,3.5) + (30:0.8)$);
\path  ($(5,3.5) + (210:0.8)$) edge[bend left =20, color = blue] ($(5,3.5) + (330:0.8)$);
\path  ($(5,3.5) + (90:0.8)$) edge[color = blue] ($(5,3.5) + (270:0.8)$);
\path  ($(7,3.5) + (-150:0.8)$) edge[bend right = 20, color = red] ($(7,3.5) + (-270:0.8)$);
\path  ($(7,3.5) + (-210:0.8)$) edge[bend right = 20, color = red] ($(7,3.5) + (-330:0.8)$);
\path  ($(7,3.5) + (-90:0.8)$) edge[bend left = 20, color = red] ($(7,3.5) + (-30:0.8)$);
\path ($(9.5,3.5) + (112.5:1.05)$)  edge[bend left =45, color = blue]  ($(9.5,3.5) + (157.5:1.05)$);
\path ($(9.5,3.5) + (135:1.05)$)  edge[bend left =20, color = blue]  ($(9.5,3.5) + (225:1.05)$);
\path ($(9.5,3.5) + (202.5:1.05)$)  edge[bend left =45, color = blue]  ($(9.5,3.5) + (247.5:1.05)$);
\path ($(9.5,3.5) + (292.5:1.05)$)  edge[bend left =20, color = red]  ($(9.5,3.5) + (45:1.05)$);
\path ($(9.5,3.5) + (315:1.05)$)  edge[bend left =90, color = red]  ($(9.5,3.5) + (337.5:1.05)$);
\path ($(9.5,3.5) + (22.5:1.05)$)  edge[bend left =45, color = red]  ($(9.5,3.5) + (67.5:1.05)$);
\end{tikzpicture}
\caption{Multiplication of chord diagrams.} \label{pict:multchord}
\end{center}
\end{figure}

The \emph{product of two chord diagrams} $C_1$ and $C_2$ is the chord diagram corresponding to the arc diagram obtained by the concatenation of two arc representations of $C_1$ and $C_2$
(see Fig.~\ref{pict:multchord})
The multiplication of chord diagrams is well-defined (i.e., the result does not depend on the arc representations) modulo the four-term relations.  

Let $C$ be a chord diagram. We denote its set of chords by $V(C)$. Let $C|_{U}$ denote the chord diagram consisting of all chords in the subset~$U\subset V(C)$. 

The comultiplication of chord diagrams is defined as
\begin{equation*}
\mu(C):= \sum_{U\subset V(C)} C|_{U}\otimes C|_{V(C)\setminus U}.\label{dias_coprod}
\end{equation*}

Multiplication and comultiplication are extended to the vector space of chord diagrams by linearity and preserve the grading.

\begin{Claim}
The multiplication and comultiplication defined above 
turn the quotient space of the vector space of chord diagrams by the four-term relations into a Hopf algebra.
\end{Claim}
We denote this Hopf algebra by $\EuScript C$. In the sequel we will call it the \emph{Hopf algebra of chord diagrams}, assuming that we consider chord diagrams up to four-term relations. 

To each chord diagram $C$ we associate its \emph{intersection graph} $\gamma(C)$.
The vertices of $\gamma(C)$ correspond to the chords of the diagram, and two vertices $v_a$ and $v_b$ in $V(\gamma(C))$ are connected by an edge if and only if the corresponding chords in $V(C)$ intersect (i.e., their endpoints $a_1, a_2, b_1,$ and $b_2$ are arranged on the circle in the order $a_1, b_1, a_2, b_2$).

\begin{Claim}[\cite{Lando}]
The map taking each chord diagram to its intersection graph extends to a graded homomorphism $\EuScript C \to \EuScript G$ of Hopf algebras. This homomorphism descends to a graded homomorphism $\EuScript C \to \EuScript F$ of Hopf algebras. 
\end{Claim}

Note that not every graph can be realized as the intersection graph of a chord diagram. Besides, two distinct chord diagrams may have the same intersection graph.

\subsection{The Hopf subalgebra generated by the complete bipartite graphs in the Hopf algebra of graphs and its Hopf subalgebras}

\begin{Def}
A graph $G$ is \emph{a complete bipartite graph} if the set $V(G)$ of its vertices can be partitioned into two subsets (\emph{parts}) $U$ and $W$ so that any two vertices $v_1 \in V(G)$ and $v_2 \in V(G)$
\begin{enumerate}
\item are not connected by an edge if they belong to the same part  {\rm(}$v_1, v_2 \in U$ or $v_1, v_2 \in W${\rm)},
\item are connected by an edge if they belong to different parts  {\rm(}$v_1 \in U$ and $v_2 \in W$ or $v_1 \in W$ and $v_2 \in U$\rm{)}.
\end{enumerate}
\end{Def}

The complete bipartite graph with parts of sizes $n$ and $m$ is denoted by $K_{n,m}$.
One of the parts of such a graph may be empty. Note that $K_{0,n}=K_{0,1}^n$.
By $\EuScript B$ we denote
the Hopf subalgebra generated by all connected complete bipartite graphs in the Hopf algebra of graphs. 
Note that any subgraph of a complete bipartite graph is also a complete bipartite graph. Thus, the vector space of complete bipartite graphs is closed under comultiplication.

By $\EuScript B^{(n)}$ we denote the Hopf subalgebra in $\EuScript B$ generated by complete bipartite graphs with no more than $n$ vertices in one of the parts, $n =0,1,2,3,\ldots$:
$$
\EuScript B^{(l)} = \left \langle K_{0,1}, K_{1,1}, K_{1,2},K_{1,3}, \ldots,K_{2,2},K_{2,3},\ldots K_{l,l}, K_{l,l+1},\ldots \right \rangle.
$$

These Hopf subalgebras are nested:
 $\EuScript B^{(0)} \subset \EuScript B^{(1)} \subset \EuScript B^{(2)} \subset \ldots\subset \EuScript B$.

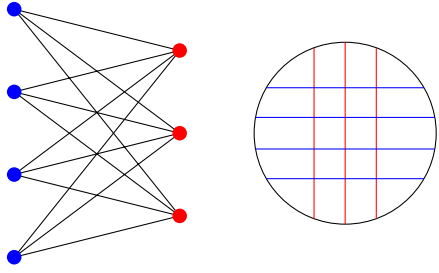
\begin{figure}[h]
\begin{center}
\begin{tikzpicture}[scale = 1.1]
\draw (1,5) -- (3,4.5) -- (1,4) -- (3,3.5) -- (1,3) -- (3,2.5) -- (1,2) -- (3,3.5) -- (1,5) -- (3, 2.5) -- (1,4);
\draw (1,2) -- (3,4.5) -- (1,3);
\fill[color = blue] (1,5) circle (2.5 pt);
\fill[color = blue] (1,4) circle (2.5 pt);
\fill[color = blue] (1,3) circle (2.5 pt);
\fill[color = blue] (1,2) circle (2.5 pt);
\fill[color = red] (3,4.5) circle (2.5 pt);
\fill[color = red] (3,3.5) circle (2.5 pt);
\fill[color = red] (3,2.5) circle (2.5 pt);

\draw  (5,3.5) circle [radius=1.1];
\path  ($(5,3.5) + (110:1.1)$) edge[color = red] ($(5,3.5) + (250:1.1)$);
\path  ($(5,3.5) + (90:1.1)$) edge[color = red] ($(5,3.5) + (270:1.1)$);
\path  ($(5,3.5) + (70:1.1)$) edge[color = red] ($(5,3.5) + (290:1.1)$);
\path  ($(5,3.5) + (150:1.1)$) edge[color = blue] ($(5,3.5) + (30:1.1)$);
\path  ($(5,3.5) + (170:1.1)$) edge[color = blue] ($(5,3.5) + (10:1.1)$);
\path  ($(5,3.5) + (190:1.1)$) edge[color = blue] ($(5,3.5) + (-10:1.1)$);
\path  ($(5,3.5) + (210:1.1)$) edge[color = blue] ($(5,3.5) + (-30:1.1)$);
\end{tikzpicture}
\end{center}
\caption{
The complete bipartite graph $K_{3,4}$ and a chord diagram with intersection graph isomorphic to  $K_{3,4}$.
\label{pic:K34}}
\end{figure}

In this paper we study the structure of the Hopf algebras $\EuScript B^{(0)}$, $\EuScript B^{(1)}$, $\EuScript B^{(2)}$, and  $\EuScript B^{(3)}$. The Hopf algebra $\EuScript B^{(0)}$ is generated by the graph  $K_{0,1}$. The Hopf algebra  $\EuScript B^{(1)}$ is generated by the graphs of the form $K_{1,n}$ (the so-called \emph{star graphs} $S_n$). The Hopf algebra $\EuScript B^{(2)}$ is generated by the graphs of the forms $K_{1,n}$ and $K_{2,n}$. The Hopf algebra $\EuScript B^{(3)}$ is generated by the graphs of the forms $K_{1,n}$,  $K_{2,n}$, and $K_{3,n}$.

\subsection{Primitive elements in the Hopf algebra of graphs}

\begin{Def}
An element  $p$ of a bialgebra is called \emph{primitive} if $\mu(p)=1\otimes p + p\otimes 1$.
\end{Def}

It is easy to show that primitive elements form a vector subspace in a bialgebra. Since any homogeneous component of a primitive element is primitive, such a vector subpace of a graded bialgebra is also graded. 

\begin{Claim}[Milnor--Moore theorem~\cite{MM}]
Over a field of characteristic zero, each connected commutative cocommutative graded bialgebra is isomorphic to the polynomial bialgebra generated by its primitive elements.
\end{Claim}

(A graded bialgebra $A_0\oplus A_1 \oplus A_2 \oplus \ldots$ over a field $\mathbb K$ is \emph{connected} if $A_0 \cong \mathbb K$.)

Decomposable elements (i.e., products of homogeneous elements of lower degree) span a vector subspace in each homogeneous subspace of a graded bialgebra. It follows from the Milnor-Moore theorem that  every homogeneous subspace is the direct sum of a subspace generated by decomposable elements and a subspace of primitive elements. Therefore, a projection~$\pi$ of each homogeneous subspace to the subspace of primitive elements along the subspace generated by decomposable elements is well defined.

\begin{Claim}[\cite{Lando, S94}] \label{cl:projG}

The projection $\pi(G)$ of any graph~$G$
on the subspace of primitive elements along the subspace generated by decomposable elements in the Hopf algebra $\EuScript G$ is given by the formula
\begin{equation}
\pi(G):=G-1! \sum_{V_1\sqcup V_2 = V(G)}G|_{V_1}\cdot G|_{V_2}+2! \sum_{V_1\sqcup V_2\sqcup V_3 = V(G)}G|_{V_1}\cdot G|_{V_2}\cdot G|_{V_3}-\ldots, \label{eq:projG}
\end{equation}
where $V(G)$ is the vertex set of $G$ and $V_1, V_2, V_3, \ldots$ are nonempty nonintersecting subsets of $V(G)$.
\end{Claim}

\subsection{
Projection onto the space of primitive elements in Hopf algebras of complete bipartite graphs.
}

The general formula \eqref{eq:projG} for the projection onto the subspace of primitive elements in the Hopf algebra of graphs is hard to use.
But this formula can be significantly simplified in some special cases. We will derive an explicit projection formula in Hopf algebras of complete bipartite graphs in terms of the generating functions
\begin{align}\label{eq:gf}
{\EuScript K}_{l}(x):=&\sum_{n=0}^{\infty} K_{l,n}\frac{x^{n+l}}{n!},\\
{\EuScript P}_{l}(x):=&\sum_{n=0}^{\infty} \pi(K_{l,n})\frac{x^{n+l}}{n!}
\notag
\end{align}
for $l = 0, 1, 2, 3$. Note that, for~$l =0$, we have
$$
{\EuScript K}_0(x)=\sum_{n=0}^{\infty} K_{0,n}\frac{x^{n}}{n!}=\sum_{n=0}^{\infty} K_{0,1}^n\frac{x^{n}}{n!}=\exp(K_{0,1}x).
$$

\begin{Remark}
Each graph $K_{l,n}$ has $n!$ automorphisms preserving the part of the graph that consists of $l$ vertices {\rm(} we call them selected vertices {\rm)}. In each summand of the generating functions  \eqref{eq:gf} the denominator $n!$ equals the number of such automorphisms of the corresponding graph, and the exponent of $x$ equals the number of vertices of the corresponding graph.
\end{Remark}

Our first main result is the following theorem.

\begin{Theorem} \label{th:projexp}

In the case of the complete bipartite graphs  $K_{0,n}$, $K_{1,n}$, $K_{2,n}$, and $K_{3,n}$, the
generating functions for projections are expressed in terms of the generating functions for graphs as follows:
\begin{align}
 \EuScript P_0(x) = &\log\EuScript K_0(x)=K_{0,1}x 
  \notag  \\ 
 \EuScript P_1(x) = &\EuScript K_1(x)\exp(- K_{0,1} x) 
  \notag \\ 
 \EuScript P_2(x) = &\EuScript K_2(x)\exp(- K_{0,1} x) - (\EuScript K_1(x)\exp(- K_{0,1} x))^2  \notag \\
				=& \EuScript K_2(x)\exp(- K_{0,1} x) - \EuScript P_1(x)^2
\label{eq:pr2_exp}\\
\EuScript P_3(x) =& \EuScript K_3(x) \exp(- K_{0,1} x) - 3\EuScript K_2(x) \EuScript K_1(x) \exp(-K_{0,1}x)^2 +2 (\EuScript K_1(x) \exp(-K_{0,1}x))^3  \notag \\
				=& \EuScript K_3(x) \exp (- K_{0,1} x) - 3\EuScript P_2(x)\EuScript P_1(x) -\EuScript P_1(x)^3
\label{eq:pr3_exp}
\end{align}
\end{Theorem}

The following definition is needed for the proof of Theorem \ref{th:projexp}.

The \emph{Stirling number of the second kind} is the number of ways to partition a set of $n$ labeled objects into $m$ nonempty unlabeled subsets. It is denoted by $\NPK{n}{m}$, $n,m \geq 0$. Here are some examples:
$$
\NPK{n}{0}=0, n \in \mathbb N;
\NPK{n}{n}=1, n\in \mathbb N \cup \lbrace 0 \rbrace;
\NPK{n}{n-1}=\binom{n}{2}, n \in \mathbb N;
\NPK{4}{2}=7.
$$

\begin{Lemma} \label{lem:altst} 
For any positive integer~$a$ and~$N$,
\begin{equation}
\sum \limits_{m=a}^{N+a} (-1)^{m-1} (m-1)! \NPK{N}{m-a} = (-1)^{a-1}(a-1)!(-a)^N. \label{eq:altst}
\end{equation}
\end{Lemma}

The proof of this lemma uses the following well-known fact (see, e.g, \cite{concrete}).
\begin{Claim}
Let $k,N \in \mathbb N \cup \lbrace 0 \rbrace$, and let $\FF{x}{k}$ denote the falling factorial:
\begin{align*}
\FF{x}{0} &:=x^0 = 1, \\
\FF{x}{k} &:= x(x-1)(x-2)\ldots (x-k+1).
\end{align*}
Then 
\begin{equation}
\sum \limits_{k=1}^{N} \NPK{N}{k} \FF{x}{k} = x^N. \label{eq:altstFF}
\end{equation}
\end{Claim}

\begin{proof}[Proof of Lemma \ref{lem:altst}]
If $a=m$, then the first summand in \eqref{eq:altst} is zero, since $\NPK{N}{0}=0$ if $N \in \mathbb N.$ Thus, we rewrite \eqref{eq:altstFF}, replacing $k$ by $m-a$, as
\begin{equation*}
\sum_{m=a}^{N+a} \NPK{N}{m-a} x(x-1)(x-2)\ldots (x-(m-a-1)) = x^N.
\end{equation*}
We substitute $x = -a$:
\begin{equation*}
\sum_{m=a}^{N+a} \NPK{N}{m-a} (-a)(-a-1)(-a-2)\ldots (-m+1) = (-a)^N.
\end{equation*}
To conclude the proof, it remains to multiply both sides by  $(-1)^{a-1}(a-1)!$.
\end{proof}

\begin{proof}[Proof of Theorem \ref{th:projexp}]

For $l=0$, there is nothing to prove. Indeed, for $n=1$, the graph~$K_{0,1}$ is a primitive element of the Hopf algebra of graphs; therefore,
$\pi(K_{0,1}) = K_{0,1}$. Further, for $n\ne1$, the graph $K_{0,n}$ is a decomposable element (since it is a disconnected graph); thus, $\pi(K_{0,n}) = 0$. 
Now consider~$l>0$.

\begin{enumerate}
\item If $l=1$, then, by Claim~\ref{cl:projG},
\begin{equation*}
\pi(K_{1,n}) = \sum_{m=1}^{n+1} (-1)^{m-1}(m-1)!\sum_{V_1\sqcup V_2 \sqcup \ldots \sqcup V_m = V(K_{1,n})}K_{1,n}|_{V_1}\cdot K_{1,n}|_{V_2} \cdot \ldots \cdot K_{1,n}|_{V_m}.
\end{equation*}
(By $V_1\sqcup V_2 \sqcup \ldots \sqcup V_m$ we denote a partition of $V(K_{1,n})$ into $m$ disjoint nonempty subsets.)

Given $m = 1,2,\ldots,n+1$, consider all possible partitions of $V(K_{1,n})$ into $m$ nonempty subsets. One of these subsets contains the selected vertex, and the remaining $m-1$ subsets do not.
We sum the products of subgraphs of $K_{1,n}$ corresponding to the partitions. For each partition, by $i$ ($0\leq i \leq n-m+1$) we denote the number of vertices belonging to the same subset as the selected vertex.
The vertices of this subset induce a subgraph  $K_{1,i}$ of  $K_{1,n}$. There are $\binom{n}{i}$ ways to choose these vertices. Furthermore, there are  $\NPK{n-i}{m-1}$ partitions of the remaining $n-i$ vertices into $m-1$ sets. There are no edges of  $K_{1,n}$ connecting these vertices. Therefore, for any such partition, the product of the corresponding graphs is just the disjoint union of   $n-i$ vertices, which is equal to the product of $n-i$ copies of $K_{0,1}$. It now follows that 
 \begin{equation}
\pi(K_{1,n}) = \sum_{m=1}^{n+1} (-1)^{m-1}(m-1)! \sum_{i =0}^{n-m+1} \binom{n}{i} \NPK{n-i}{m-1} K_{1,i}K_{0,1}^{n - i}.\label{eq:proj1}
\end{equation}

If $i>n-m+1$, then $\NPK{n-i}{m-1} = 0$, so we can change the summation limits in the second sum (without changing the value of \eqref{eq:proj1}) as \label{item:lims}
\begin{equation*}
\pi(K_{1,n}) = \sum_{m=1}^{n+1} (-1)^{m-1}(m-1)! \sum_{i =0}^{n} \binom{n}{i} \NPK{n-i}{m-1} K_{1,i}K_{0,1}^{n - i}.\label{eq:proj1lim=n}
\end{equation*}

Changing the order of summation, we obtain
\begin{equation*}
\pi(K_{1,n}) =\sum_{i =0}^{n} K_{1,i}K_{0,1}^{n - i}\binom{n}{i} \sum_{m=1}^{n+1} (-1)^{m-1}(m-1)!   \NPK{n-i}{m-1} .\label{eq:proj1ordChanged}
\end{equation*}

The application of Lemma~\ref{lem:altst} for $a = 1$ yields
\begin{equation}
\pi(K_{1,n}) =\sum_{i =0}^{n} K_{1,i}K_{0,1}^{n - i}\binom{n}{i}  (-1)^{n-i} .\label{eq:proj1(-1)}
\end{equation}

On the other hand,
\begin{equation*}
\EuScript K_1(x)\exp(-K_{0,1} x) = \left(\sum_{n=0}^{\infty} K_{1,n}\frac{x^{n+1}}{n!}\right ) \left( \sum_{j=0}^{\infty} \frac{(-K_{0,1} x)^j}{j!}\right ).
\end{equation*}
The coefficient of $x^{n+1}$ in this expression equals
\begin{equation*}
\sum_{i=0}^n K_{1,i}\frac{(-K_{0,1})^{n-i}}{(n-i)!i!}.
\end{equation*}
Therefore, the coefficient of $\frac{x^{n+1}}{n!}$ equals \eqref{eq:proj1(-1)}, as required.

\item Suppose that $l =2$.
Now we consider the projections of $K_{2,n}$ and partitions of $V(K_{2,n})$. There are two cases: in the first case, both selected vertices belong to the same part, and in the second, they belong to two different parts. Thus, ${\EuScript P}_2$ is the sum of the two expressions corresponding to these cases. For the first case, the reasoning is the same as  for $K_{1,n}$. Therefore, the first summand in ${\EuScript P}_2$ equals $\EuScript K_2(x)\exp(-K_{0,1} x)$.
Given a graph $K_{2,n}$ and $m = 1, \ldots, n+1$,
the second summand equals
\begin{equation*}
	\sum_{m=2}^{n+2} (-1)^{m-1}(m-1)! \sum_{k=0}^n \sum_{i=0}^k \binom{n}{i,k-i,n-k}\NPK{n-k}{m-2}K_{1,i} K_{1,k-i} K_{0,1}^{n-k}, \label{eq:proj2}
\end{equation*}
where $k$ is the total number of vertices in the two parts containing selected vertices, and $i$ and $k-i$ are the numbers of vertices in each of these parts.

Now we argue as in case~\ref{item:lims}. We change the summation limits, replacing $0\leq k \leq n-m+2$ by $0\leq k \leq n$, and the order of summation. Then, in view of the fact that $\NPK{n-k}{m-2} = 0$ for $m>n-k$, we again change the summation limits, replacing $2 \leq m \leq n+2$ by $2 \leq m \leq n- k +2$. We obtain

\begin{equation*}
 \sum_{k=0}^n \sum_{i=0}^k \binom{n}{i,k-i,n-k} K_{0,1}^{n-k}K_{1,i}K_{1,k-i} \sum_{m=2}^{n - k+2}(-1)^{m-1}(m-1)!\NPK{n-k}{m-2}.
\end{equation*}

The application of Lemma~\ref{lem:altst} for $a = 2$ yields
\begin{equation}
-\sum_{k=0}^n \sum_{i=0}^k  \binom{n}{i,k-i,n-k} K_{1,i}K_{1,k-i}K_{0,1}^{n-k} (-2)^{n-k}. \label{eq:projAlt2}
\end{equation}

On the other hand, let us consider the coefficient of $\frac{x^{n+2}}{n!}$ in the second summand on the right-hand side of \eqref{eq:pr2_exp}. This summand equals
\begin{multline*}
- (\EuScript K_1(x)\exp(-K_{0,1} x))^2 =
 - (\EuScript K_1(x)^2 \exp(-2 K_{0,1} x)) = \\
-   \sum_{i=0}^{\infty} K_{1,i} \frac{x^{i+1}}{i!} \sum_{j=0}^{\infty} K_{1,j} \frac{x^{j+1}}{j!} \sum_{a=0}^{\infty} \frac{(-2)^{a} K_{0,1}^a x^a}{a!},
\end{multline*}
and the coefficient of $\frac{x^{n+2}}{n!}$ in this expression equals
\begin{equation*}
- n! \sum_{k=0}^n \sum_{i=0}^k K_{1,i} K_{1,k-i} \frac{1}{i!(k-i)!(n-k)!} (-2)^{n-k} K_{0,1}^{n-k}.
\end{equation*}
This is equal to \eqref{eq:projAlt2}, which concludes the proof of $\eqref{eq:pr2_exp}$.

 \item Let $l=3$.  The projection of $K_{3,n}$ equals

\begin{multline*}
\pi(K_{3,n}) = \sum_{m=1}^{n+1} (-1)^{m-1} (m-1)! \sum_{k=0}^{n-m+1} \binom{n}{k}\NPK{n-k}{m-1}K_{3,k}K_{0,1}^{n-k} \label{eq:proj3} \\
+ \sum_{m=2}^{n+2} (-1)^{m-1}(m-1)! \binom{3}{2} \sum_{k=0}^{n-m+2} \sum_{i=0}^k \binom{n}{i,k-i,n-k} \NPK{n-k}{m-2}K_{2,i}K_{1,k-i}K_{0,1}^{n-k} \\
+ \sum_{m=3}^{n+3} (-1)^{m-1}(m-1)! \sum_{k=0}^{n-m+3} \sum_{i=0}^k \sum_{j=0}^{k-i} \binom{n}{i,j,k-(i+j),n-k} \NPK{n-k}{m-3}K_{1,i}K_{1,j}K_{1,k-(i+j)}K_{0,1}^{n-k}.
\end{multline*}
Here the second summand corresponds to the partitions of $V(K_{3,n})$ such that two of the selected vertices belong to one part and the third, to another part. The binomial coefficient $\binom{3}{2}$  in the second summand is the number of ways to choose these two vertices.
The further argument is the same as for $l = 1, 2$.
\end{enumerate}
\end{proof}

\section{The $\sltwo$ Weight System}\label{sec:sl2}
In this section we compute the values of the $\sltwo$ weight system at the chord diagrams 
with intersection graph $K_{l,n}$, $l \leq 3$.
Then we use Theorem~\ref{th:projexp} to compute the values 
of the $\sltwo$ weight system at the projections of such complete bipartite graphs on the space of primitive elements.

The results which we obtain confirm a conjecture due to Lando, which states that
if $G$ is the intersection graph of a chord diagram such that the length of the longest cycle of $G$ (the circumference of $G$) 
is at most $2 l, l\ge1$,
then  the value of the  $\sltwo$ weight system at the projection of this chord diagram on the subspace of primitive elements
is a polynomial of degree at most~$l$. 

\subsection{Definition of the $\sltwo$ weight system}

\label{subseq:defth}

Let $R$ be a ring,and let $A$ be an algebra over $R$. A linear function $w \colon \EuScript C \to A$ vanishing at any four-term element is called a \emph{weight system} on $\EuScript C$.
We consider only the case where
$R = {\mathbb C}$ and $A={\mathbb C}[c]$.

Let $\mathfrak g$ be a Lie algebra of finite dimension over $\mathbb C$ endowed with a nondegenerate bilinear invariant  form  $(\cdot , \cdot)$. (A form is invariant if $(\lbrack x ,y \rbrack, z) = (x ,\lbrack y,z \rbrack)$ for any $x,y,z \in \mathfrak g$.) Let $X = \lbrace x_1, x_2, \ldots, x_m \rbrace$ be an orthonormal basis of $\mathfrak g$ with respect to this form. We use $U(\mathfrak g)$ to denote the universal enveloping algebra of the Lie algebra $\mathfrak g$.
Consider the map $w_{\mathfrak g} \colon \EuScript C \to U(\mathfrak g)$ defined as follows.

Suppose given a chord diagram $C$ and an arc representation $a$ of $C$. Let $V(a)$ be the set of all arcs of $a$, and let $\nu$ be a map $\nu \colon V(a) \to \lbrace 1, 2, \ldots, m \rbrace$. 
With the diagram $a$ and the map $\nu$ we associate the element $w_X(a,\nu) \in U(\mathfrak g)$ obtained as follows: at both ends of each arc $v\in V(a)$ we write the element $x_{\nu(v)} \in X$ and multiply all written elements from left to right. We denote this product by $w_X(a, \nu)$ and the sum of such products over all possible maps by $w_X(a)$:
\begin{equation}
w_X(a) := \sum \limits_{\nu} w_X(a,\nu). \label{wsldef}
\end{equation}

\begin{figure}
\begin{center}
\begin{tikzpicture}
\draw  (2,3.5) circle [radius=1.2];
\draw (4,2.5)--(14,2.5);
\node(center) at (2,3.5) {};
\node(origin) at (4,2.5) {};
\path  ($(center) + (30:1.2)$) edge[bend left =20, color = black] ($(center) + (150:1.2)$);
\path  ($(origin) + (30*0.027,0)$) edge[bend left =90, color = black] ($(origin) + (150*0.027,0)$);
\node at ($(origin) + (30*0.027,-0.25)$) {$x_{i_1}$};
\node at ($(origin) + (160*0.027,-0.25)$) {$x_{i_1}$};
\path  ($(center) + (60:1.2)$) edge[bend left =20, color = teal] ($(center) + (120:1.2)$);
\path  ($(origin) + (60*0.027,0)$) edge[bend left =90, color = teal] ($(origin) + (120*0.027,0)$);
\node at ($(origin) + (60*0.027,-0.25)$) {$x_{i_2}$};
\node at ($(origin) + (118*0.027,-0.25)$) {$x_{i_2}$};
\path  ($(center) + (90:1.2)$) edge[bend left =20, color = olive] ($(center) + (240:1.2)$);
\path  ($(origin) + (90*0.027,0)$) edge[bend left =90, color = olive] ($(origin) + (240*0.027,0)$);
\node at ($(origin) + (90*0.027,-0.25)$) {$x_{i_3}$};
\node at ($(origin) + (240*0.027,-0.25)$) {$x_{i_3}$};
\path  ($(center) + (135:1.2)$) edge[bend left =20, color = purple] ($(center) + (270:1.2)$);
\path  ($(origin) + (135*0.027,0)$) edge[bend left =90, color = purple] ($(origin) + (270*0.027,0)$);
\node at ($(origin) + (141*0.027,-0.25)$) {$x_{i_4}$};
\node at ($(origin) + (270*0.027,-0.25)$) {$x_{i_4}$};
\path  ($(center) + (210:1.2)$) edge[bend left =20, color = green] ($(center) + (330:1.2)$);
\path  ($(origin) + (210*0.027,0)$) edge[bend left =90, color = green] ($(origin) + (330*0.027,0)$);
\node at ($(origin) + (210*0.027,-0.25)$) {$x_{i_5}$};
\node at ($(origin) + (330*0.027,-0.25)$) {$x_{i_5}$};
\end{tikzpicture}
\caption{
\label{pic:arc-sl2}
Computation of the value of the weight system that corresponds to a Lie algebra with orthonormal basis $x_1,\dots,x_m$ at an arc representation of a chord diagram.
}
\end{center}
\end{figure}
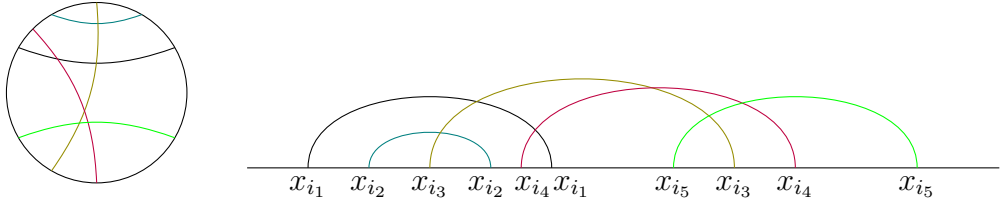

For example, the value of the weight system corresponding to a Lie algebra with orthonormal basis $x_1,\dots, x_m$ at the arc diagram shown in Fig.~\ref{pic:arc-sl2} equals
$$
\sum_{i_1=1}^m\sum_{i_2=1}^m
\sum_{i_3=1}^m
\sum_{i_4=1}^m
\sum_{i_5=1}^m
x_{i_1}x_{i_2}x_{i_3}x_{i_2}x_{i_4}x_{i_1}x_{i_5}x_{i_3}x_{i_4}x_{i_5}.
$$

\begin{Claim}
\begin{enumerate}
\item For any $C \in \EuScript C$ the result of this operation is determined uniquely and does not depend on the choice of an arc representation of $C$.
\item For any  $a$, $w_X(a) \in Z(U(\mathfrak g))$, where $Z(U(\mathfrak g))$ is the center of the universal enveloping algebra.
\item The element $w_X(a)$ does not depend on the choice of an orthonormal basis.
\item The map from chord diagrams to $ Z(U(\mathfrak g))$ thus defined satisfies the four-term relations. Therefore, it extends to a homomorphism of commutative algebras. 
\end{enumerate}
\end{Claim}

Since we define the multiplication of chord diagrams as the concatenation of its arc representations, the weight system corresponding to any Lie algebra is multiplicative.

We consider this construction for the simplest case of a noncommutative Lie algebra, namely, for the Lie algebra $\sltwo$.
This Lie algebra is generated by three elements $x,y,$ and $z$ with relations
\begin{gather*}
\lbrack x, y \rbrack = z, \\
\lbrack  y, z \rbrack = x,\\
\lbrack z, x \rbrack = y.
\end{gather*}
The bilinear form is given by the relations
\begin{gather*}
(x,x) = (y,y) = (z,z) = 1,\\
(x,y) = (y,z) = (z,x) = 0.
\end{gather*}
The center of the universal enveloping algebra $Z(U(\sltwo))$ is isomorphic to the algebra of polynomials in the Casimir element $c=x^2 + y^2 + z^2\in U(\sltwo)$. Hence the formula  \eqref{wsldef} defines a map $w_{\sltwo} \colon \EuScript C \to \mathbb C\lbrack c \rbrack$. This map is a homomorphism of algebras and is called the \emph{$\sltwo$ weight system on $\EuScript C$.}

This definition immediately implies the following assertion.
\begin{conseq}
The value of the weight system  $\sltwo$ on a chord diagram with only one chord equals $c$.
\end{conseq}

In~\cite{ChL} the following nontrivial statement, which links the $\sltwo$ weight system with polynomial graph invariants, was proved. Note that its analogue for more complicated Lie algebras, e.g, for $\mathfrak{sl}_3$, turns out to be wrong.

\begin{Claim}
The value of the $\sltwo$ weight system at a chord diagram depends only on the intersection graph of this diagram.
\end{Claim}

To compute the values of the $\sltwo$ weight system at chord diagrams, we use the multiplicativity of this weight system
and the \emph{Chmutov--Varchenko recurrence relations}.
For simplicity, we will identify the value of the
 $\sltwo$ weight system at a chord diagram with the diagram itself. 
As in Fig.~\ref{pict:4term}, diagrams may contain other chords with endpoints on the dashed arcs; the sets of these additional chords must be the same for all terms of each equation.

\begin{Claim}[Chmutov--Varchenko reccurence relations, \cite{bib:ChV}] \label{cl:ChV}
Let $D$ be a chord diagram. Assume that its intersection graph is connected, i.e, $D$ is not a product of two chord diagrams of lower order. Then the following assertions hold.
\begin{enumerate}
\item If $D$ contains a \emph{leaf}, i.e., a chord intersecting only one chord, then
\begin{equation}
w_{\sltwo}(D) = (c-1)w_{\sltwo}(D'), \label{eq:ChV_leaf}
\end{equation}
where $D'$ is the chord diagram obtained from $D$ by deleting the leaf:
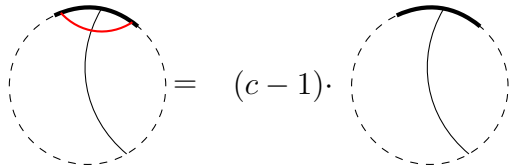
\begin{center}
\begin{tikzpicture}[scale = 1.3]
\draw[dashed]  (1,3.5) circle [radius=0.8];
\path  ($(1,3.5) + (115:0.8)$) edge[ultra thick, bend left =32] ($(1,3.5) + (50:0.8)$);
\node  at (2,3.5) {$=$};
\node  at (3,3.5) {$(c-1)\cdot$};
\draw[dashed]  (4.5,3.5) circle [radius=0.8];
\node at (5.5,3.25) {.};
\path  ($(4.5,3.5) + (115:0.8)$) edge[ultra thick, bend left =32] ($(4.5,3.5) + (50:0.8)$);
\path  ($(1,3.5) + (80:0.8)$) edge[bend right =40] ($(1,3.5) + (300:0.8)$);
\path  ($(1,3.5) + (110:0.8)$) edge[thick, color = red, bend right =40] ($(1,3.5) + (55:0.8)$);
\path  ($(4.5,3.5) + (80:0.8)$) edge[bend right =40] ($(4.5,3.5) + (300:0.8)$);
\end{tikzpicture}
\end{center}

\item
The following equations hold:

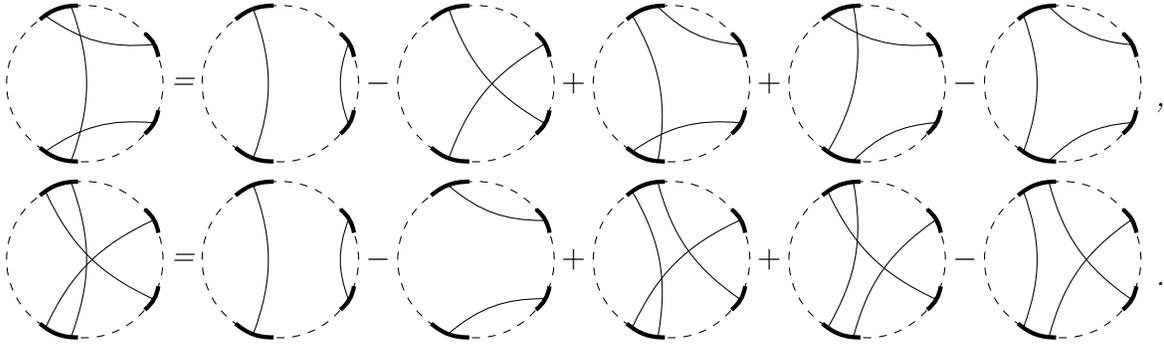
\begin{tikzpicture}[scale = 1.3]
\draw[dashed]  (1,3.5) circle [radius=0.8];
\node  at (2,3.5) {=};
\draw[dashed]  (3,3.5) circle [radius=0.8];
\node  at (4,3.5) {$-$};
\draw[dashed]  (5,3.5) circle [radius=0.8];
\node  at (6,3.5) {$+$};
\draw[dashed]  (7,3.5) circle [radius=0.8];
\node  at (8,3.5) {$+$};
\draw[dashed]  (9,3.5) circle [radius=0.8];
\node  at (10,3.5) {$-$};
\draw[dashed]  (11,3.5) circle [radius=0.8];
\node at (12,3.25) {,};

\foreach \n in {0,2,...,10}{
\node(center) at (1 + \n,3.5) {};
\path ($(center) + (20:0.8)$) edge[bend right =20, ultra thick] ($(center) + (40:0.8)$);
\path ($(center) + (95:0.8)$) edge[bend right =20, ultra thick] ($(center) + (125:0.8)$);
\path ($(center) + (-20:0.8)$) edge[bend right =-20, ultra thick] ($(center) + (-40:0.8)$);
\path ($(center) + (-95:0.8)$) edge[bend right =-20, ultra thick] ($(center) + (-125:0.8)$);
}
\path  ($(1,3.5) + (120:0.8)$) edge[bend right =20] ($(1,3.5) + (30:0.8)$);
\path  ($(1,3.5) + (-120:0.8)$) edge[bend left =20] ($(1,3.5) + (330:0.8)$);
\path  ($(1,3.5) + (100:0.8)$) edge[bend left = 20] ($(1,3.5) + (260:0.8)$);
\path  ($(3,3.5) + (110:0.8)$) edge[bend left =20] ($(3,3.5) + (-110:0.8)$);
\path  ($(3,3.5) + (30:0.8)$) edge[bend right =20] ($(3,3.5) + (330:0.8)$);
\path  ($(5,3.5) + (110:0.8)$) edge[bend right =20]  ($(5,3.5) + (330:0.8)$);
\path  ($(5,3.5) + (30:0.8)$) edge[bend right =20] ($(5,3.5) + (-110:0.8)$);
\path  ($(7,3.5) + (120:0.8)$) edge[bend left = 20] ($(7,3.5) + (260:0.8)$);
\path  ($(7,3.5) + (-120:0.8)$) edge[bend left = 20] ($(7,3.5) + (330:0.8)$);
\path  ($(7,3.5) + (100:0.8)$) edge[bend right = 20] ($(7,3.5) + (30:0.8)$);
\path  ($(9,3.5) + (-120:0.8)$) edge[bend right = 20] ($(9,3.5) + (-260:0.8)$);
\path  ($(9,3.5) + (120:0.8)$) edge[bend right = 20] ($(9,3.5) + (-330:0.8)$);
\path  ($(9,3.5) + (-100:0.8)$) edge[bend left = 20] ($(9,3.5) + (-30:0.8)$);
\path  ($(11,3.5) + (30:0.8)$) edge[bend left = 20] ($(11,3.5) + (100:0.8)$);
\path  ($(11,3.5) + (120:0.8)$) edge[bend left = 20] ($(11,3.5) + (240:0.8)$);
\path  ($(11,3.5) + (260:0.8)$) edge[bend left = 20] ($(11,3.5) + (330:0.8)$);
\end{tikzpicture}


\begin{tikzpicture}[scale = 1.3]
\draw[dashed]  (1,3.5) circle [radius=0.8];
\node  at (2,3.5) {=};
\draw[dashed]  (3,3.5) circle [radius=0.8];
\node  at (4,3.5) {$-$};
\draw[dashed]  (5,3.5) circle [radius=0.8];
\node  at (6,3.5) {$+$};
\draw[dashed]  (7,3.5) circle [radius=0.8];
\node  at (8,3.5) {$+$};
\draw[dashed]  (9,3.5) circle [radius=0.8];
\node  at (10,3.5) {$-$};
\draw[dashed]  (11,3.5) circle [radius=0.8];
\node at (12,3.25) {.};


\foreach \n in {0,2,...,10}{
\node(center) at (1 + \n,3.5) {};
\path ($(center) + (20:0.8)$) edge[bend right =20, ultra thick] ($(center) + (40:0.8)$);
\path ($(center) + (95:0.8)$) edge[bend right =20, ultra thick] ($(center) + (125:0.8)$);
\path ($(center) + (-20:0.8)$) edge[bend right =-20, ultra thick] ($(center) + (-40:0.8)$);
\path ($(center) + (-95:0.8)$) edge[bend right =-20, ultra thick] ($(center) + (-125:0.8)$);
}
\path  ($(1,3.5) + (120:0.8)$) edge[bend right =20] ($(1,3.5) + (330:0.8)$);
\path  ($(1,3.5) + (240:0.8)$) edge[bend left =20] ($(1,3.5) + (30:0.8)$);
\path  ($(1,3.5) + (100:0.8)$) edge[bend left = 20] ($(1,3.5) + (260:0.8)$);
\path  ($(3,3.5) + (110:0.8)$) edge[bend left =20] ($(3,3.5) + (-110:0.8)$);
\path  ($(3,3.5) + (330:0.8)$) edge[bend left =20] ($(3,3.5) + (30:0.8)$);
\path  ($(5,3.5) + (110:0.8)$) edge[bend right =20]  ($(5,3.5) + (30:0.8)$);
\path  ($(5,3.5) + (330:0.8)$) edge[bend right =20] ($(5,3.5) + (-110:0.8)$);
\path  ($(7,3.5) + (120:0.8)$) edge[bend left = 20] ($(7,3.5) + (260:0.8)$);
\path  ($(7,3.5) + (-120:0.8)$) edge[bend left = 20] ($(7,3.5) + (30:0.8)$);
\path  ($(7,3.5) + (100:0.8)$) edge[bend right = 20] ($(7,3.5) + (330:0.8)$);
\path  ($(9,3.5) + (-120:0.8)$) edge[bend right = 20] ($(9,3.5) + (-260:0.8)$);
\path  ($(9,3.5) + (120:0.8)$) edge[bend right = 20] ($(9,3.5) + (-30:0.8)$);
\path  ($(9,3.5) + (-100:0.8)$) edge[bend left = 20] ($(9,3.5) + (-330:0.8)$);
\path  ($(11,3.5) + (330:0.8)$) edge[bend left = 20] ($(11,3.5) + (100:0.8)$);
\path  ($(11,3.5) + (120:0.8)$) edge[bend left = 20] ($(11,3.5) + (240:0.8)$);
\path  ($(11,3.5) + (260:0.8)$) edge[bend left = 20] ($(11,3.5) + (30:0.8)$);
\end{tikzpicture}

\end{enumerate}
\end{Claim}

Now we will compute the values of the $\sltwo$ weight system on complete bipartite graphs with restricted size of one of the parts, using these reccurence relations.

\subsection{The values of the $\sltwo$ weight system at the graphs $K_{1,n}$, $K_{2,n}$, and $K_{3,n}$.}

Relation~\eqref{eq:ChV_leaf} implies the following assertion
\begin{conseq}
If the intersection graph of a chord diagram is a tree with $n$ vertices, then the value of the $\sltwo$ weight system at this chord diagram equals  $c(c-1)^{n-1}.$
\end{conseq}

In particular, $w_\sltwo(K_{1,n})=c(c-1)^n$ for  $n=0,1,2,\dots$.

The next theorem is the our main result about values of the $\sltwo$ weight system. Since its value at a chord diagram depends only on the intersection graph of this chord diagram, we use graphs as the arguments of the weight system. We set $k_{i,j} = w_\sltwo(K_{i,j})$, $i,j=0,1,2,3,\dots$.

\begin{Theorem}
\label{th:rec}
The following relations hold:
\begin{align}
k_{0,n} =& c^n, 
 \notag \\
k_{1,n} =& c(c-1)^n, \label{eq:W1} \\
k_{2,n} =& (c-3)k_{2,n-1} +c (c^n + (c-1)^{n-1}) {\text{, }} n \geq 2, \label{eq:n2rec} \\
k_{3,n} =& - 2 (c+3)k_{3,n-1} + 3c(c-6)k_{3,n-2} \notag  \\
+&2k_{2,n+2} + 11k_{2,n+1} - (12c^2 - 11c - 16)k_{2,n}  \notag \\ 
+& (16c^3-55c^2+32c+9)k_{2,n-1} -3c(2c^3-11c^2+16c-9)k_{2,n-2} \notag \\
+& c(c-1)^{n-2}(4c-1)(3c^2-2c+1) - 8c^{n+1} { \text{, }} n\geq 3.  \label{eq:n3rec}
\end{align}
\end{Theorem}

\begin{proof} [Proof of the Theorem~\ref{th:rec}]
\begin{enumerate}
\item The equation $k_{0,n} = c^n$  follows from the multiplicativity of the weight system.
\item The equation $k_{1,n} = c(c-1)^n$ follows from the fact that the complete bipartite graph $K_{1,n}$ is a tree.
\item

To prove the recurrence relation~ \eqref{eq:n2rec} for~$k_{2,n}$, we use the multiplicativity of the weight system and the Chmutov--Varchenko relations (see Claim \ref{cl:ChV}).

Consider the following chord diagram:

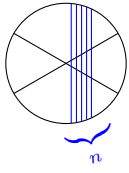
\begin{tikzpicture}[scale = 1]
\draw  (5,0) circle [radius=0.8];
\node(center) at (1 + 4,0) {};
\path  ($(5,0) + (150:0.8)$) edge ($(5,0) + (330:0.8)$);
\path  ($(5,0) + (30:0.8)$) edge ($(5,0) + (210:0.8)$);
\path  ($(1+ 4,0) + (85:0.8)$) edge[color = blue] ($(1 + 4,0) + (275:0.8)$);
\path  ($(1 + 4,0) + (80:0.8)$) edge[color = blue] ($(1 + 4,0) + (280:0.8)$);
\path  ($(1 + 4,0) + (75:0.8)$) edge[color = blue] ($(1 + 4,0) + (285:0.8)$);
\path  ($(1 + 4,0) + (70:0.8)$) edge[color = blue] ($(1 + 4,0) + (290:0.8)$);
\path  ($(1 + 4,0) + (65:0.8)$) edge[color = blue] ($(1 + 4,0) + (295:0.8)$);
\node[color = blue, rotate = 17] at (1.35 + 4,-1.1) {\footnotesize{$\underbrace{\phantom{}}_{n}$}};
\end{tikzpicture}

The intersection graph of this chord diagram is a complete tripartite graph with parts of sizes $1,1,$ and $n$ and hence we denote this chord diagram by $K_{1,1,n}$ and the value of the $\sltwo$ weight system at it by $k_{1,1,n}$.

The first of the Chmutov--Varchenko six-term relations for a chord diagram with intersection graph~$K_{2,n}$ is

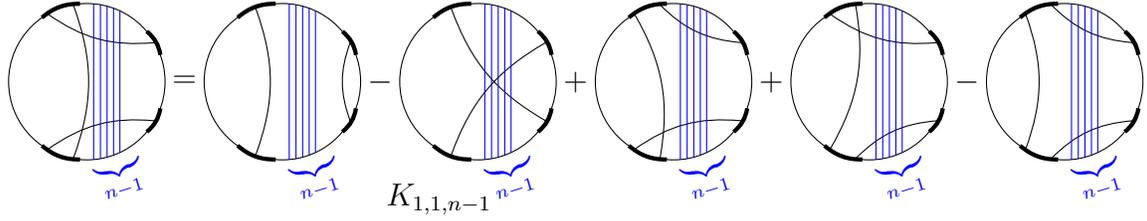
\begin{tikzpicture} [scale = 1.3]
\draw  (1,3.5) circle [radius=0.8];
\node  at (2,3.5) {=};
\draw  (3,3.5) circle [radius=0.8];
\node  at (4,3.5) {$-$};
\draw  (5,3.5) circle [radius=0.8];
\node  at (6,3.5) {$+$};
\draw  (7,3.5) circle [radius=0.8];
\node  at (8,3.5) {$+$};
\draw  (9,3.5) circle [radius=0.8];
\node  at (10,3.5) {$-$};
\draw  (11,3.5) circle [radius=0.8];

\foreach \n in {0,2,...,10}{
\node(center) at (1 + \n,3.5) {};
\path ($(center) + (20:0.8)$) edge[bend right =20, ultra thick] ($(center) + (40:0.8)$);
\path ($(center) + (95:0.8)$) edge[bend right =20, ultra thick] ($(center) + (125:0.8)$);
\path ($(center) + (-20:0.8)$) edge[bend right =-20, ultra thick] ($(center) + (-40:0.8)$);
\path ($(center) + (-95:0.8)$) edge[bend right =-20, ultra thick] ($(center) + (-125:0.8)$);
}
\path  ($(1,3.5) + (120:0.8)$) edge[bend right =20] ($(1,3.5) + (30:0.8)$);
\path  ($(1,3.5) + (-120:0.8)$) edge[bend left =20] ($(1,3.5) + (330:0.8)$);
\path  ($(1,3.5) + (100:0.8)$) edge[bend left = 20] ($(1,3.5) + (260:0.8)$);
\path  ($(3,3.5) + (110:0.8)$) edge[bend left =20] ($(3,3.5) + (-110:0.8)$);
\path  ($(3,3.5) + (30:0.8)$) edge[bend right =20] ($(3,3.5) + (330:0.8)$);
\path  ($(5,3.5) + (110:0.8)$) edge[bend right =20]  ($(5,3.5) + (330:0.8)$);
\path  ($(5,3.5) + (30:0.8)$) edge[bend right =20] ($(5,3.5) + (-110:0.8)$);
\path  ($(7,3.5) + (120:0.8)$) edge[bend left = 20] ($(7,3.5) + (260:0.8)$);
\path  ($(7,3.5) + (-120:0.8)$) edge[bend left = 20] ($(7,3.5) + (330:0.8)$);
\path  ($(7,3.5) + (100:0.8)$) edge[bend right = 20] ($(7,3.5) + (30:0.8)$);
\path  ($(9,3.5) + (-120:0.8)$) edge[bend right = 20] ($(9,3.5) + (-260:0.8)$);
\path  ($(9,3.5) + (120:0.8)$) edge[bend right = 20] ($(9,3.5) + (-330:0.8)$);
\path  ($(9,3.5) + (-100:0.8)$) edge[bend left = 20] ($(9,3.5) + (-30:0.8)$);
\path  ($(11,3.5) + (30:0.8)$) edge[bend left = 20] ($(11,3.5) + (100:0.8)$);
\path  ($(11,3.5) + (120:0.8)$) edge[bend left = 20] ($(11,3.5) + (240:0.8)$);
\path  ($(11,3.5) + (260:0.8)$) edge[bend left = 20] ($(11,3.5) + (330:0.8)$);
\foreach \n in {0,2,...,10}{
\path  ($(1+ \n,3.5) + (85:0.8)$) edge[color = blue] ($(1 + \n,3.5) + (275:0.8)$);
\path  ($(1 + \n,3.5) + (80:0.8)$) edge[color = blue] ($(1 + \n,3.5) + (280:0.8)$);
\path  ($(1 + \n,3.5) + (75:0.8)$) edge[color = blue] ($(1 + \n,3.5) + (285:0.8)$);
\path  ($(1 + \n,3.5) + (70:0.8)$) edge[color = blue] ($(1 + \n,3.5) + (290:0.8)$);
\path  ($(1 + \n,3.5) + (65:0.8)$) edge[color = blue] ($(1 + \n,3.5) + (295:0.8)$);
\node[color = blue, rotate = 17] at (1.35 + \n,2.5) {\footnotesize{$\underbrace{\phantom{}}_{n-1}$}};
}
\node at (4.6, 2.3) {$K_{1,1,n-1}$};
\end{tikzpicture}

The sum of the last three summands equals $(c-2) k_{2,n-1}$. The first summand on the right-hand side equals $c^{n+1}$. Therefore,
\begin{equation}
k_{2,n} = c^{n+1} - k_{1,1,n-1} + (c-2) k_{2,n-1}.
\label{eq:Xn}
\end{equation}
Now we will use the following Chmutov--Varchenko relation to compute the value of $k_{1,1,n}$:

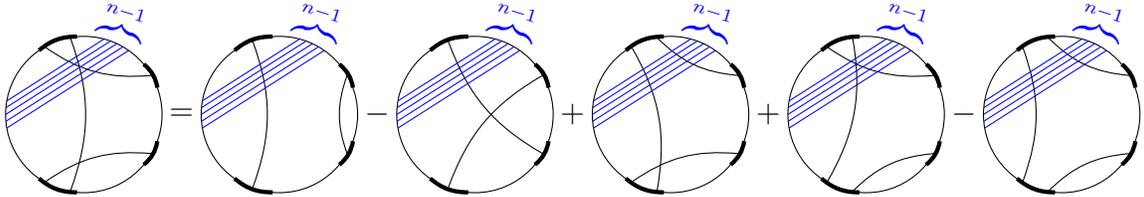
\begin{tikzpicture}[scale = 1.3]
\draw  (1,3.5) circle [radius=0.8];
\node  at (2,3.5) {=};
\draw  (3,3.5) circle [radius=0.8];
\node  at (4,3.5) {$-$};
\draw  (5,3.5) circle [radius=0.8];
\node  at (6,3.5) {$+$};
\draw  (7,3.5) circle [radius=0.8];
\node  at (8,3.5) {$+$};
\draw  (9,3.5) circle [radius=0.8];
\node  at (10,3.5) {$-$};
\draw  (11,3.5) circle [radius=0.8];

\foreach \n in {0,2,...,10}{
\node(center) at (1 + \n,3.5) {};
\path ($(center) + (20:0.8)$) edge[bend right =20, ultra thick] ($(center) + (40:0.8)$);
\path ($(center) + (95:0.8)$) edge[bend right =20, ultra thick] ($(center) + (125:0.8)$);
\path ($(center) + (-20:0.8)$) edge[bend right =-20, ultra thick] ($(center) + (-40:0.8)$);
\path ($(center) + (-95:0.8)$) edge[bend right =-20, ultra thick] ($(center) + (-125:0.8)$);
}
\path  ($(1,3.5) + (120:0.8)$) edge[bend right =20] ($(1,3.5) + (30:0.8)$);
\path  ($(1,3.5) + (-120:0.8)$) edge[bend left =20] ($(1,3.5) + (330:0.8)$);
\path  ($(1,3.5) + (100:0.8)$) edge[bend left = 20] ($(1,3.5) + (260:0.8)$);
\path  ($(3,3.5) + (110:0.8)$) edge[bend left =20] ($(3,3.5) + (-110:0.8)$);
\path  ($(3,3.5) + (30:0.8)$) edge[bend right =20] ($(3,3.5) + (330:0.8)$);
\path  ($(5,3.5) + (110:0.8)$) edge[bend right =20]  ($(5,3.5) + (330:0.8)$);
\path  ($(5,3.5) + (30:0.8)$) edge[bend right =20] ($(5,3.5) + (-110:0.8)$);
\path  ($(7,3.5) + (120:0.8)$) edge[bend left = 20] ($(7,3.5) + (260:0.8)$);
\path  ($(7,3.5) + (-120:0.8)$) edge[bend left = 20] ($(7,3.5) + (330:0.8)$);
\path  ($(7,3.5) + (100:0.8)$) edge[bend right = 20] ($(7,3.5) + (30:0.8)$);
\path  ($(9,3.5) + (-120:0.8)$) edge[bend right = 20] ($(9,3.5) + (-260:0.8)$);
\path  ($(9,3.5) + (120:0.8)$) edge[bend right = 20] ($(9,3.5) + (-330:0.8)$);
\path  ($(9,3.5) + (-100:0.8)$) edge[bend left = 20] ($(9,3.5) + (-30:0.8)$);
\path  ($(11,3.5) + (30:0.8)$) edge[bend left = 20] ($(11,3.5) + (100:0.8)$);
\path  ($(11,3.5) + (120:0.8)$) edge[bend left = 20] ($(11,3.5) + (240:0.8)$);
\path  ($(11,3.5) + (260:0.8)$) edge[bend left = 20] ($(11,3.5) + (330:0.8)$);
\foreach \n in {0,2,...,10}{
\path  ($(1+ \n,3.5) + (170:0.8)$) edge[color = blue] ($(1 + \n,3.5) + (75:0.8)$);
\path  ($(1 + \n,3.5) + (175:0.8)$) edge[color = blue] ($(1 + \n,3.5) + (70:0.8)$);
\path  ($(1 + \n,3.5) + (180:0.8)$) edge[color = blue] ($(1 + \n,3.5) + (65:0.8)$);
\path  ($(1 + \n,3.5) + (185:0.8)$) edge[color = blue] ($(1 + \n,3.5) + (60:0.8)$);
\path  ($(1 + \n,3.5) + (190:0.8)$) edge[color = blue] ($(1 + \n,3.5) + (55:0.8)$);
\node[color = blue, rotate = -17] at (1.4 + \n,4.45) {\footnotesize{$\overbrace{\phantom{}}^{n-1}$}};
}
\end{tikzpicture}

We obtain the following relation between $k_{1,1,n-1}$ and $k_{2,n-1}$:

\begin{equation}
(c-1)k_{1,1,n-1} = ck_{1,n-1} - k_{1,n} + (c-1) k_{2,n-1} + c k_{1,1,n-1} - c k_{2,n-1};
\end{equation}
therefore,
$$
k_{1,1,n-1}=k_{2,n-1}-c(c-1)^{n-1}.
$$
To conclude the proof of~\ref{eq:n2rec}, it remains to substitute this expression for $k_{1,1,n}$ into~\eqref{eq:Xn}.

\item
Now let us prove~\eqref{eq:n3rec}. The following equations follow from the Chmutov--Varchenko recurrence relations:
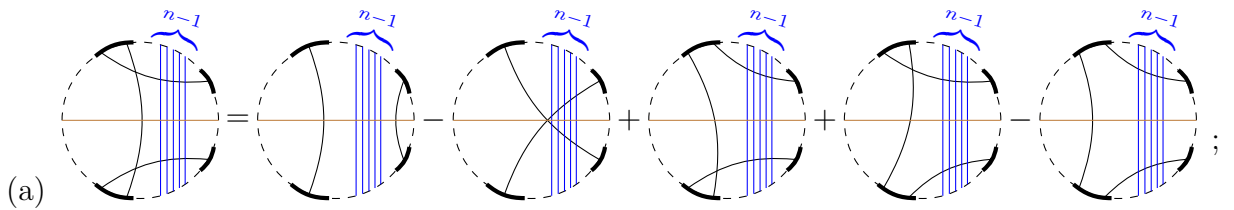
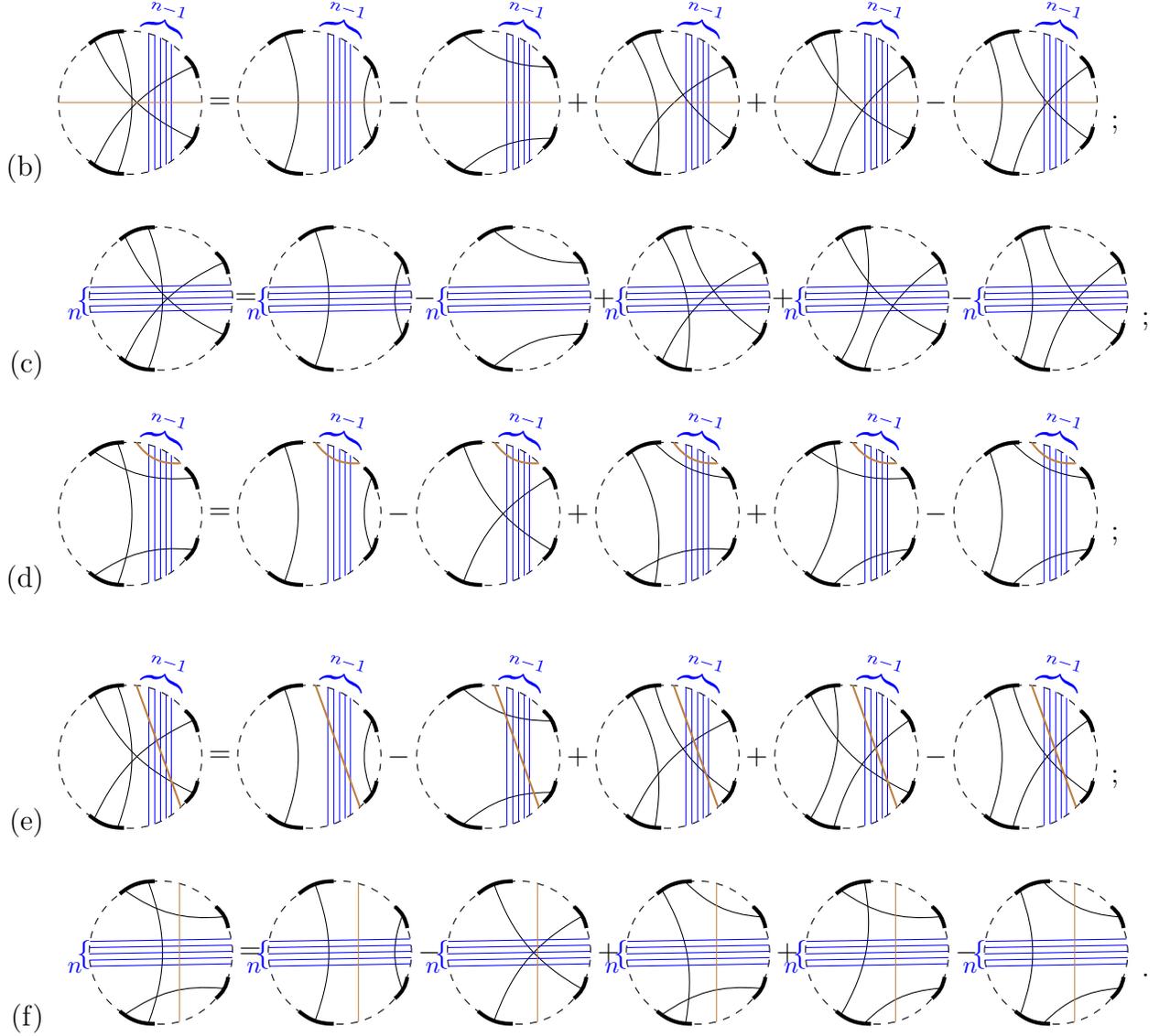
\begin{enumerate}
\item \begin{tikzpicture}[scale = 1.3]
\draw[dashed]  (1,3.5) circle [radius=0.8];
\node  at (2,3.5) {=};
\draw[dashed]  (3,3.5) circle [radius=0.8];
\node  at (4,3.5) {$-$};
\draw[dashed]  (5,3.5) circle [radius=0.8];
\node  at (6,3.5) {$+$};
\draw[dashed]  (7,3.5) circle [radius=0.8];
\node  at (8,3.5) {$+$};
\draw[dashed]  (9,3.5) circle [radius=0.8];
\node  at (10,3.5) {$-$};
\draw[dashed]  (11,3.5) circle [radius=0.8];
\node at (12,3.25) {;};

\foreach \n in {0,2,...,10}{
\node(center) at (1 + \n,3.5) {};
\path ($(center) + (20:0.8)$) edge[bend right =20, ultra thick] ($(center) + (40:0.8)$);
\path ($(center) + (95:0.8)$) edge[bend right =20, ultra thick] ($(center) + (125:0.8)$);
\path ($(center) + (-20:0.8)$) edge[bend right =-20, ultra thick] ($(center) + (-40:0.8)$);
\path ($(center) + (-95:0.8)$) edge[bend right =-20, ultra thick] ($(center) + (-125:0.8)$);
}

\path  ($(1,3.5) + (120:0.8)$) edge[bend right =20] ($(1,3.5) + (30:0.8)$);
\path  ($(1,3.5) + (-120:0.8)$) edge[bend left =20] ($(1,3.5) + (330:0.8)$);
\path  ($(1,3.5) + (100:0.8)$) edge[bend left = 20] ($(1,3.5) + (260:0.8)$);
\path  ($(3,3.5) + (110:0.8)$) edge[bend left =20] ($(3,3.5) + (-110:0.8)$);
\path  ($(3,3.5) + (30:0.8)$) edge[bend right =20] ($(3,3.5) + (330:0.8)$);
\path  ($(5,3.5) + (110:0.8)$) edge[bend right =20]  ($(5,3.5) + (330:0.8)$);
\path  ($(5,3.5) + (30:0.8)$) edge[bend right =20] ($(5,3.5) + (-110:0.8)$);
\path  ($(7,3.5) + (120:0.8)$) edge[bend left = 20] ($(7,3.5) + (260:0.8)$);
\path  ($(7,3.5) + (-120:0.8)$) edge[bend left = 20] ($(7,3.5) + (330:0.8)$);
\path  ($(7,3.5) + (100:0.8)$) edge[bend right = 20] ($(7,3.5) + (30:0.8)$);
\path  ($(9,3.5) + (-120:0.8)$) edge[bend right = 20] ($(9,3.5) + (-260:0.8)$);
\path  ($(9,3.5) + (120:0.8)$) edge[bend right = 20] ($(9,3.5) + (-330:0.8)$);
\path  ($(9,3.5) + (-100:0.8)$) edge[bend left = 20] ($(9,3.5) + (-30:0.8)$);
\path  ($(11,3.5) + (30:0.8)$) edge[bend left = 20] ($(11,3.5) + (100:0.8)$);
\path  ($(11,3.5) + (120:0.8)$) edge[bend left = 20] ($(11,3.5) + (240:0.8)$);
\path  ($(11,3.5) + (260:0.8)$) edge[bend left = 20] ($(11,3.5) + (330:0.8)$);
\foreach \n in {0,2,...,10}{
\path  ($(1+ \n,3.5) + (75:0.8)$) edge[color = blue] ($(1 + \n,3.5) + (-75:0.8)$);
\path  ($(1 + \n,3.5) + (70:0.8)$) edge[color = blue] ($(1 + \n,3.5) + (-70:0.8)$);
\path  ($(1 + \n,3.5) + (65:0.8)$) edge[color = blue] ($(1 + \n,3.5) + (-65:0.8)$);
\path  ($(1 + \n,3.5) + (60:0.8)$) edge[color = blue] ($(1 + \n,3.5) + (-60:0.8)$);
\path  ($(1 + \n,3.5) + (55:0.8)$) edge[color = blue] ($(1 + \n,3.5) + (-55:0.8)$);
\node[color = blue, rotate = -17] at (1.4 + \n,4.45) {\footnotesize{$\overbrace{\phantom{}}^{n-1}$}};
\path  ($(1 + \n,3.5) + (180:0.8)$) edge[color = brown] ($(1 + \n,3.5) + (0:0.8)$);
}
\end{tikzpicture}
 \label{item:Eq31}
\item \begin{tikzpicture}[scale = 1.3]
\draw[dashed]  (1,3.5) circle [radius=0.8];
\node  at (2,3.5) {=};
\draw[dashed]  (3,3.5) circle [radius=0.8];
\node  at (4,3.5) {$-$};
\draw[dashed]  (5,3.5) circle [radius=0.8];
\node  at (6,3.5) {$+$};
\draw[dashed]  (7,3.5) circle [radius=0.8];
\node  at (8,3.5) {$+$};
\draw[dashed]  (9,3.5) circle [radius=0.8];
\node  at (10,3.5) {$-$};
\draw[dashed]  (11,3.5) circle [radius=0.8];
\node at (12,3.25) {;};

\foreach \n in {0,2,...,10}{
\node(center) at (1 + \n,3.5) {};
\path ($(center) + (20:0.8)$) edge[bend right =20, ultra thick] ($(center) + (40:0.8)$);
\path ($(center) + (95:0.8)$) edge[bend right =20, ultra thick] ($(center) + (125:0.8)$);
\path ($(center) + (-20:0.8)$) edge[bend right =-20, ultra thick] ($(center) + (-40:0.8)$);
\path ($(center) + (-95:0.8)$) edge[bend right =-20, ultra thick] ($(center) + (-125:0.8)$);
}

\path  ($(1,3.5) + (120:0.8)$) edge[bend right =20] ($(1,3.5) + (330:0.8)$);
\path  ($(1,3.5) + (240:0.8)$) edge[bend left =20] ($(1,3.5) + (30:0.8)$);
\path  ($(1,3.5) + (100:0.8)$) edge[bend left = 20] ($(1,3.5) + (260:0.8)$);
\path  ($(3,3.5) + (110:0.8)$) edge[bend left =20] ($(3,3.5) + (-110:0.8)$);
\path  ($(3,3.5) + (330:0.8)$) edge[bend left =20] ($(3,3.5) + (30:0.8)$);
\path  ($(5,3.5) + (110:0.8)$) edge[bend right =20]  ($(5,3.5) + (30:0.8)$);
\path  ($(5,3.5) + (330:0.8)$) edge[bend right =20] ($(5,3.5) + (-110:0.8)$);
\path  ($(7,3.5) + (120:0.8)$) edge[bend left = 20] ($(7,3.5) + (260:0.8)$);
\path  ($(7,3.5) + (-120:0.8)$) edge[bend left = 20] ($(7,3.5) + (30:0.8)$);
\path  ($(7,3.5) + (100:0.8)$) edge[bend right = 20] ($(7,3.5) + (330:0.8)$);
\path  ($(9,3.5) + (-120:0.8)$) edge[bend right = 20] ($(9,3.5) + (-260:0.8)$);
\path  ($(9,3.5) + (120:0.8)$) edge[bend right = 20] ($(9,3.5) + (-30:0.8)$);
\path  ($(9,3.5) + (-100:0.8)$) edge[bend left = 20] ($(9,3.5) + (-330:0.8)$);
\path  ($(11,3.5) + (330:0.8)$) edge[bend left = 20] ($(11,3.5) + (100:0.8)$);
\path  ($(11,3.5) + (120:0.8)$) edge[bend left = 20] ($(11,3.5) + (240:0.8)$);
\path  ($(11,3.5) + (260:0.8)$) edge[bend left = 20] ($(11,3.5) + (30:0.8)$);
\foreach \n in {0,2,...,10}{
\path  ($(1+ \n,3.5) + (75:0.8)$) edge[color = blue] ($(1 + \n,3.5) + (-75:0.8)$);
\path  ($(1 + \n,3.5) + (70:0.8)$) edge[color = blue] ($(1 + \n,3.5) + (-70:0.8)$);
\path  ($(1 + \n,3.5) + (65:0.8)$) edge[color = blue] ($(1 + \n,3.5) + (-65:0.8)$);
\path  ($(1 + \n,3.5) + (60:0.8)$) edge[color = blue] ($(1 + \n,3.5) + (-60:0.8)$);
\path  ($(1 + \n,3.5) + (55:0.8)$) edge[color = blue] ($(1 + \n,3.5) + (-55:0.8)$);
\node[color = blue, rotate = -17] at (1.4 + \n,4.45) {\footnotesize{$\overbrace{\phantom{}}^{n-1}$}};
\path  ($(1 + \n,3.5) + (180:0.8)$) edge[color = brown] ($(1 + \n,3.5) + (0:0.8)$);
}
\end{tikzpicture}
\label{item:Eq32}
\item  \begin{tikzpicture}[scale = 1.3]
\draw[dashed]  (1,3.5) circle [radius=0.8];
\node  at (1.95,3.5) {=};
\draw[dashed]  (3,3.5) circle [radius=0.8];
\node  at (3.95,3.5) {$-$};
\draw[dashed]  (5,3.5) circle [radius=0.8];
\node  at (5.95,3.5) {$+$};
\draw[dashed]  (7,3.5) circle [radius=0.8];
\node  at (7.95,3.5) {$+$};
\draw[dashed]  (9,3.5) circle [radius=0.8];
\node  at (9.95,3.5) {$-$};
\draw[dashed]  (11,3.5) circle [radius=0.8];
\node at (12,3.25) {;};

\foreach \n in {0,2,...,10}{
\node(center) at (1 + \n,3.5) {};
\path ($(center) + (20:0.8)$) edge[bend right =20, ultra thick] ($(center) + (40:0.8)$);
\path ($(center) + (95:0.8)$) edge[bend right =20, ultra thick] ($(center) + (125:0.8)$);
\path ($(center) + (-20:0.8)$) edge[bend right =-20, ultra thick] ($(center) + (-40:0.8)$);
\path ($(center) + (-95:0.8)$) edge[bend right =-20, ultra thick] ($(center) + (-125:0.8)$);
}

\path  ($(1,3.5) + (120:0.8)$) edge[bend right =20] ($(1,3.5) + (330:0.8)$);
\path  ($(1,3.5) + (240:0.8)$) edge[bend left =20] ($(1,3.5) + (30:0.8)$);
\path  ($(1,3.5) + (100:0.8)$) edge[bend left = 20] ($(1,3.5) + (260:0.8)$);
\path  ($(3,3.5) + (110:0.8)$) edge[bend left =20] ($(3,3.5) + (-110:0.8)$);
\path  ($(3,3.5) + (330:0.8)$) edge[bend left =20] ($(3,3.5) + (30:0.8)$);
\path  ($(5,3.5) + (110:0.8)$) edge[bend right =20]  ($(5,3.5) + (30:0.8)$);
\path  ($(5,3.5) + (330:0.8)$) edge[bend right =20] ($(5,3.5) + (-110:0.8)$);
\path  ($(7,3.5) + (120:0.8)$) edge[bend left = 20] ($(7,3.5) + (260:0.8)$);
\path  ($(7,3.5) + (-120:0.8)$) edge[bend left = 20] ($(7,3.5) + (30:0.8)$);
\path  ($(7,3.5) + (100:0.8)$) edge[bend right = 20] ($(7,3.5) + (330:0.8)$);
\path  ($(9,3.5) + (-120:0.8)$) edge[bend right = 20] ($(9,3.5) + (-260:0.8)$);
\path  ($(9,3.5) + (120:0.8)$) edge[bend right = 20] ($(9,3.5) + (-30:0.8)$);
\path  ($(9,3.5) + (-100:0.8)$) edge[bend left = 20] ($(9,3.5) + (-330:0.8)$);
\path  ($(11,3.5) + (330:0.8)$) edge[bend left = 20] ($(11,3.5) + (100:0.8)$);
\path  ($(11,3.5) + (120:0.8)$) edge[bend left = 20] ($(11,3.5) + (240:0.8)$);
\path  ($(11,3.5) + (260:0.8)$) edge[bend left = 20] ($(11,3.5) + (30:0.8)$);
\foreach \n in {0,2,...,10}{
\path  ($(1+ \n,3.5) + (171:0.8)$) edge[color = blue] ($(1 + \n,3.5) + (11:0.8)$);
\path  ($(1 + \n,3.5) + (176:0.8)$) edge[color = blue] ($(1 + \n,3.5) + (6:0.8)$);
\path  ($(1 + \n,3.5) + (181:0.8)$) edge[color = blue] ($(1 + \n,3.5) + (1:0.8)$);
\path  ($(1 + \n,3.5) + (186:0.8)$) edge[color = blue] ($(1 + \n,3.5) + (-4:0.8)$);
\path  ($(1 + \n,3.5) + (191:0.8)$) edge[color = blue] ($(1 + \n,3.5) + (-9:0.8)$);
\node[color = blue] at (0.14+\n,3.475) {$\mathbf{\lbrace }$};
\node[color = blue] at (0.04+\n,3.33) {\footnotesize{$n$}};
}
\end{tikzpicture}

\label{item:Eq33}
\item \begin{tikzpicture}[scale = 1.3]
\draw[dashed]  (1,3.5) circle [radius=0.8];
\node  at (2,3.5) {=};
\draw[dashed]  (3,3.5) circle [radius=0.8];
\node  at (4,3.5) {$-$};
\draw[dashed]  (5,3.5) circle [radius=0.8];
\node  at (6,3.5) {$+$};
\draw[dashed]  (7,3.5) circle [radius=0.8];
\node  at (8,3.5) {$+$};
\draw[dashed]  (9,3.5) circle [radius=0.8];
\node  at (10,3.5) {$-$};
\draw[dashed]  (11,3.5) circle [radius=0.8];
\node at (12,3.25) {;};

\foreach \n in {0,2,...,10}{
\node(center) at (1 + \n,3.5) {};
\path ($(center) + (20:0.8)$) edge[bend right =20, ultra thick] ($(center) + (40:0.8)$);
\path ($(center) + (95:0.8)$) edge[bend right =20, ultra thick] ($(center) + (125:0.8)$);
\path ($(center) + (-20:0.8)$) edge[bend right =-20, ultra thick] ($(center) + (-40:0.8)$);
\path ($(center) + (-95:0.8)$) edge[bend right =-20, ultra thick] ($(center) + (-125:0.8)$);
}

\path  ($(1,3.5) + (120:0.8)$) edge[bend right =20] ($(1,3.5) + (30:0.8)$);
\path  ($(1,3.5) + (-120:0.8)$) edge[bend left =20] ($(1,3.5) + (330:0.8)$);
\path  ($(1,3.5) + (100:0.8)$) edge[bend left = 20] ($(1,3.5) + (260:0.8)$);
\path  ($(3,3.5) + (110:0.8)$) edge[bend left =20] ($(3,3.5) + (-110:0.8)$);
\path  ($(3,3.5) + (30:0.8)$) edge[bend right =20] ($(3,3.5) + (330:0.8)$);
\path  ($(5,3.5) + (110:0.8)$) edge[bend right =20]  ($(5,3.5) + (330:0.8)$);
\path  ($(5,3.5) + (30:0.8)$) edge[bend right =20] ($(5,3.5) + (-110:0.8)$);
\path  ($(7,3.5) + (120:0.8)$) edge[bend left = 20] ($(7,3.5) + (260:0.8)$);
\path  ($(7,3.5) + (-120:0.8)$) edge[bend left = 20] ($(7,3.5) + (330:0.8)$);
\path  ($(7,3.5) + (100:0.8)$) edge[bend right = 20] ($(7,3.5) + (30:0.8)$);
\path  ($(9,3.5) + (-120:0.8)$) edge[bend right = 20] ($(9,3.5) + (-260:0.8)$);
\path  ($(9,3.5) + (120:0.8)$) edge[bend right = 20] ($(9,3.5) + (-330:0.8)$);
\path  ($(9,3.5) + (-100:0.8)$) edge[bend left = 20] ($(9,3.5) + (-30:0.8)$);
\path  ($(11,3.5) + (30:0.8)$) edge[bend left = 20] ($(11,3.5) + (100:0.8)$);
\path  ($(11,3.5) + (120:0.8)$) edge[bend left = 20] ($(11,3.5) + (240:0.8)$);
\path  ($(11,3.5) + (260:0.8)$) edge[bend left = 20] ($(11,3.5) + (330:0.8)$);
\foreach \n in {0,2,...,10}{
\path  ($(1+ \n,3.5) + (75:0.8)$) edge[color = blue] ($(1 + \n,3.5) + (-75:0.8)$);
\path  ($(1 + \n,3.5) + (70:0.8)$) edge[color = blue] ($(1 + \n,3.5) + (-70:0.8)$);
\path  ($(1 + \n,3.5) + (65:0.8)$) edge[color = blue] ($(1 + \n,3.5) + (-65:0.8)$);
\path  ($(1 + \n,3.5) + (60:0.8)$) edge[color = blue] ($(1 + \n,3.5) + (-60:0.8)$);
\path  ($(1 + \n,3.5) + (55:0.8)$) edge[color = blue] ($(1 + \n,3.5) + (-55:0.8)$);
\node[color = blue, rotate = -17] at (1.4 + \n,4.45) {\footnotesize{$\overbrace{\phantom{}}^{n-1}$}};
\path  ($(1 + \n,3.5) + (85:0.8)$) edge[bend right = 30, thick, color = brown] ($(1 + \n,3.5) + (45:0.8)$);
}
\end{tikzpicture}
\label{item:Eq34a}
\item  \begin{tikzpicture}[scale = 1.3]
\draw[dashed]  (1,3.5) circle [radius=0.8];
\node  at (2,3.5) {=};
\draw[dashed]  (3,3.5) circle [radius=0.8];
\node  at (4,3.5) {$-$};
\draw[dashed]  (5,3.5) circle [radius=0.8];
\node  at (6,3.5) {$+$};
\draw[dashed]  (7,3.5) circle [radius=0.8];
\node  at (8,3.5) {$+$};
\draw[dashed]  (9,3.5) circle [radius=0.8];
\node  at (10,3.5) {$-$};
\draw[dashed]  (11,3.5) circle [radius=0.8];
\node at (12,3.25) {;};

\foreach \n in {0,2,...,10}{
\node(center) at (1 + \n,3.5) {};
\path ($(center) + (20:0.8)$) edge[bend right =20, ultra thick] ($(center) + (40:0.8)$);
\path ($(center) + (95:0.8)$) edge[bend right =20, ultra thick] ($(center) + (125:0.8)$);
\path ($(center) + (-20:0.8)$) edge[bend right =-20, ultra thick] ($(center) + (-40:0.8)$);
\path ($(center) + (-95:0.8)$) edge[bend right =-20, ultra thick] ($(center) + (-125:0.8)$);
}

\path  ($(1,3.5) + (120:0.8)$) edge[bend right =20] ($(1,3.5) + (330:0.8)$);
\path  ($(1,3.5) + (240:0.8)$) edge[bend left =20] ($(1,3.5) + (30:0.8)$);
\path  ($(1,3.5) + (100:0.8)$) edge[bend left = 20] ($(1,3.5) + (260:0.8)$);
\path  ($(3,3.5) + (110:0.8)$) edge[bend left =20] ($(3,3.5) + (-110:0.8)$);
\path  ($(3,3.5) + (330:0.8)$) edge[bend left =20] ($(3,3.5) + (30:0.8)$);
\path  ($(5,3.5) + (110:0.8)$) edge[bend right =20]  ($(5,3.5) + (30:0.8)$);
\path  ($(5,3.5) + (330:0.8)$) edge[bend right =20] ($(5,3.5) + (-110:0.8)$);
\path  ($(7,3.5) + (120:0.8)$) edge[bend left = 20] ($(7,3.5) + (260:0.8)$);
\path  ($(7,3.5) + (-120:0.8)$) edge[bend left = 20] ($(7,3.5) + (30:0.8)$);
\path  ($(7,3.5) + (100:0.8)$) edge[bend right = 20] ($(7,3.5) + (330:0.8)$);
\path  ($(9,3.5) + (-120:0.8)$) edge[bend right = 20] ($(9,3.5) + (-260:0.8)$);
\path  ($(9,3.5) + (120:0.8)$) edge[bend right = 20] ($(9,3.5) + (-30:0.8)$);
\path  ($(9,3.5) + (-100:0.8)$) edge[bend left = 20] ($(9,3.5) + (-330:0.8)$);
\path  ($(11,3.5) + (330:0.8)$) edge[bend left = 20] ($(11,3.5) + (100:0.8)$);
\path  ($(11,3.5) + (120:0.8)$) edge[bend left = 20] ($(11,3.5) + (240:0.8)$);
\path  ($(11,3.5) + (260:0.8)$) edge[bend left = 20] ($(11,3.5) + (30:0.8)$);
\foreach \n in {0,2,...,10}{
\path  ($(1+ \n,3.5) + (75:0.8)$) edge[color = blue] ($(1 + \n,3.5) + (-75:0.8)$);
\path  ($(1 + \n,3.5) + (70:0.8)$) edge[color = blue] ($(1 + \n,3.5) + (-70:0.8)$);
\path  ($(1 + \n,3.5) + (65:0.8)$) edge[color = blue] ($(1 + \n,3.5) + (-65:0.8)$);
\path  ($(1 + \n,3.5) + (60:0.8)$) edge[color = blue] ($(1 + \n,3.5) + (-60:0.8)$);
\path  ($(1 + \n,3.5) + (55:0.8)$) edge[color = blue] ($(1 + \n,3.5) + (-55:0.8)$);
\node[color = blue, rotate = -17] at (1.4 + \n,4.45) {\footnotesize{$\overbrace{\phantom{}}^{n-1}$}};
\path  ($(1 + \n,3.5) + (85:0.8)$) edge[thick, color = brown] ($(1 + \n,3.5) + (-45:0.8)$);
}
\end{tikzpicture}
\label{item:Eq35a}
\item \begin{tikzpicture}[scale = 1.3]
\draw[dashed]  (1,3.5) circle [radius=0.8];
\node  at (2,3.5) {=};
\draw[dashed]  (3,3.5) circle [radius=0.8];
\node  at (4,3.5) {$-$};
\draw[dashed]  (5,3.5) circle [radius=0.8];
\node  at (6,3.5) {$+$};
\draw[dashed]  (7,3.5) circle [radius=0.8];
\node  at (8,3.5) {$+$};
\draw[dashed]  (9,3.5) circle [radius=0.8];
\node  at (10,3.5) {$-$};
\draw[dashed]  (11,3.5) circle [radius=0.8];
\node at (12,3.25) {.};

\foreach \n in {0,2,...,10}{
\node(center) at (1 + \n,3.5) {};
\path ($(center) + (20:0.8)$) edge[bend right =20, ultra thick] ($(center) + (40:0.8)$);
\path ($(center) + (95:0.8)$) edge[bend right =20, ultra thick] ($(center) + (125:0.8)$);
\path ($(center) + (-20:0.8)$) edge[bend right =-20, ultra thick] ($(center) + (-40:0.8)$);
\path ($(center) + (-95:0.8)$) edge[bend right =-20, ultra thick] ($(center) + (-125:0.8)$);
}

\path  ($(1,3.5) + (120:0.8)$) edge[bend right =20] ($(1,3.5) + (30:0.8)$);
\path  ($(1,3.5) + (-120:0.8)$) edge[bend left =20] ($(1,3.5) + (330:0.8)$);
\path  ($(1,3.5) + (100:0.8)$) edge[bend left = 20] ($(1,3.5) + (260:0.8)$);
\path  ($(3,3.5) + (110:0.8)$) edge[bend left =20] ($(3,3.5) + (-110:0.8)$);
\path  ($(3,3.5) + (30:0.8)$) edge[bend right =20] ($(3,3.5) + (330:0.8)$);
\path  ($(5,3.5) + (110:0.8)$) edge[bend right =20]  ($(5,3.5) + (330:0.8)$);
\path  ($(5,3.5) + (30:0.8)$) edge[bend right =20] ($(5,3.5) + (-110:0.8)$);
\path  ($(7,3.5) + (120:0.8)$) edge[bend left = 20] ($(7,3.5) + (260:0.8)$);
\path  ($(7,3.5) + (-120:0.8)$) edge[bend left = 20] ($(7,3.5) + (330:0.8)$);
\path  ($(7,3.5) + (100:0.8)$) edge[bend right = 20] ($(7,3.5) + (30:0.8)$);
\path  ($(9,3.5) + (-120:0.8)$) edge[bend right = 20] ($(9,3.5) + (-260:0.8)$);
\path  ($(9,3.5) + (120:0.8)$) edge[bend right = 20] ($(9,3.5) + (-330:0.8)$);
\path  ($(9,3.5) + (-100:0.8)$) edge[bend left = 20] ($(9,3.5) + (-30:0.8)$);
\path  ($(11,3.5) + (30:0.8)$) edge[bend left = 20] ($(11,3.5) + (100:0.8)$);
\path  ($(11,3.5) + (120:0.8)$) edge[bend left = 20] ($(11,3.5) + (240:0.8)$);
\path  ($(11,3.5) + (260:0.8)$) edge[bend left = 20] ($(11,3.5) + (330:0.8)$);
\foreach \n in {0,2,...,10}{
\path  ($(1+ \n,3.5) + (171:0.8)$) edge[color = blue] ($(1 + \n,3.5) + (11:0.8)$);
\path  ($(1 + \n,3.5) + (176:0.8)$) edge[color = blue] ($(1 + \n,3.5) + (6:0.8)$);
\path  ($(1 + \n,3.5) + (181:0.8)$) edge[color = blue] ($(1 + \n,3.5) + (1:0.8)$);
\path  ($(1 + \n,3.5) + (186:0.8)$) edge[color = blue] ($(1 + \n,3.5) + (-4:0.8)$);
\path  ($(1 + \n,3.5) + (191:0.8)$) edge[color = blue] ($(1 + \n,3.5) + (-9:0.8)$);
\node[color = blue] at (0.14+\n,3.475) {$\mathbf{\lbrace }$};
\node[color = blue] at (0.04+\n,3.33) {\footnotesize{$n$}};
\path  ($(1 + \n,3.5) + (75:0.8)$) edge[color = brown] ($(1 + \n,3.5) + (-75:0.8)$);
}
\end{tikzpicture}

\label{item:Eq36}
\end{enumerate}
We have six linear equations. Each chord diagram in these equations contains a collection of $n$ nonintersecting chords represented by the parallel line segments in the figure. If we add Eqs.~\ref{item:Eq31}--\ref{item:Eq36}  with $n-1$ instead of $n$ chords in each such collection and Eqs. \ref{item:Eq33} and \ref{item:Eq36} with $n-2$ chords in each such collection, we obtain a system of $14$ linear equations. 
We eliminate 13 variables and obtain the recurrence relation~\eqref{eq:n3rec}. 
\end{enumerate}
\end{proof}

Solving the recurrence relations \eqref{eq:n2rec} and \eqref{eq:n3rec}, we obtain the values of the $\sltwo$ weight system at the complete bipartite graphs $K_{2,n}$ and $K_{3,n}$.

\begin{conseq} 
The following relations hold:
\begin{align}
k_{2,n} &= \frac{1}{6} c ((4c-3) (c-3)^n + 3 (c-1)^n+ 2 c^{n+1})  \label{eq:W2} \\
k_{3,n} &= \frac{1}{30} c \bigl (3(4c^2-11c+6) (c-6)^n \notag \\
&\phantom{=\frac{1}{30} c \bigl (} +10(4c- 3) (c-3)^n + 6(3c^2-2c+2) (c-1)^n+ 5 c^{n+1}\bigr).
 \label{eq:W3}
\end{align}
\end{conseq}

Let $\ti{\EuScript K}_{l}$ denote the exponential generating functions for the values of the weight system $w_{\sltwo}$ at the complete bipartite graphs  $K_{{l},n}$ $(l = 1,2,3)$.

Relations \eqref{eq:W1}, \eqref{eq:W2}, and \eqref{eq:W3} give the following expressions for these generating functions.
\begin{conseq} \label{conseq:WKGen} 
The following relations hold:
\begin{align}
\ti {\EuScript K}_1(x)& = x c e^{(c-1)x}  
\notag \\
\ti {\EuScript K}_2(x) &=x^2 \frac{c}{6}\left ( (4c-3)e^{(c-3)x}+3e^{(c-1)x} +2ce^{cx} \right)  
\notag \\
\ti {\EuScript K}_3(x) &= x^3\frac{c}{30} \left( 3(4 c^2 - 11 c + 6) e^{(c - 6) x}  + 10(4 c - 3) e^{(c - 3) x} + 6(3 c^2 - 2 c + 2) e^{(c - 1) x} + 5ce^{cx}\right). \notag 
\end{align}
\end{conseq}

\subsection{The values of the $\sltwo$ weight system at the projections $\pi(K_{l,n})$, $l = 1,2,3$.}

Let $\ti {\EuScript P}_l$ denote the exponential generating function for the values of the $\sltwo$ weight system at the projections of the complete bipartite graphs $K_{l,n}$ on the subspace of primitive elements.
Combining Theorem~\ref{th:projexp} and Corollary~\ref{conseq:WKGen}, we obtain expressions for $\ti{\EuScript P}_{l}$, $l = 1,2,3$.
\begin{conseq}
The following relations hold:
\begin{align}
\ti {\EuScript P}_1(x) =& cx e^{-x}, \label{eq:WProj1n}\\
\ti {\EuScript P}_2(x) =& \frac16cx^2 \left((4c-3)e^{-3x}-6ce^{-2x}+ 3 e^{-x} + 2c  \right),
\label{eq:WProj2n}
\\
\ti {\EuScript P}_3(x) =& \frac{1}{30} cx^3 ((12 c^2 - 33 c + 18) e^{-6 x}
  - c (60 c - 45) e^{-4 x} + (60 c^2 + 40 c - 30) e^{-3 x} \notag \\
 -& 45 ce^{-2 x} - (12 c^2 + 12 c - 12) e^{-x} + 5 c). \label{eq:WProj3n}
\end{align}
\end{conseq}

For each $\ti {\EuScript P}_l(x)$, the coefficient of each exponent is a polynomial in $c$ of degree at most $l$. This implies the following statement.

\begin{Claim}
The following assertions hold:
\begin{enumerate}
\item The value $w_{\sltwo}(\pi(K_{1,n}))$ equals $(-1)^nc$
\item  The value $w_{\sltwo}(\pi(K_{2,n}))$ is a polynomial of degree less than or equal to $2$. If $n \geq 2$, then the degree of this polynomial equals 2.
\item The value $w_{\sltwo}(\pi(K_{3,n}))$ is a polynomial of degree less than or equal to $3$. If $n \geq 3$, then the degree of this polynomial equals 3.
\end{enumerate}
\label{cl:w12}
\end{Claim}

This proves the following conjecture in the particular case of the graphs $K_{1,n}$, $K_{2,n}$, and $K_{3,n}$.

\begin{Conjecture}[S.K.Lando]
Suppose given a chord diagram $C$. Let $\pi(C)$ be its projection on the space of primitive elements, and let $\gamma(C)$ be its intersection graph. Then $w_{\sltwo}(\pi(C))$ is a polynomial of degree less than or equal to one half of the circumference (i.e., the length of the longest cycle) of $\gamma(C)$.
\end{Conjecture}

For the case of complete bipartite graphs, this conjecture can be restated as follows.
\begin{Conjecture}[S.K.~Lando]
The value of the $\sltwo$ weight system at the complete bipartite graph $K_{l,n}$ is a polynomial of degree $\min(l,n)$.
\end{Conjecture}

We also give formulas for the ordinary generating functions for the values of $w_{\sltwo}$ at the projections of complete bipartite graphs onto the subspace of primitive elements. We denote such a generating function by  $\ti{P}_l$; then
$$
\ti{P}_l(s)=\sum_{n=0}^\infty w_{\sltwo}(\pi(K_{l,n}))s^n.
$$

Using the fact that the ordinary generating function corresponding to the exponential generating fuction  $e^{ax}$ is $\frac1{1-as}$, we obtain the following result.

\begin{Theorem}
The following relations hold:
\begin{align*}
\ti{P}_1(s)&=\frac{cs}{1+s},\\
\ti{P}_2(s)&=\frac{cs((1+2s)+2cs(1+s))}{(1+s)(1+2s)(1+3s)},\\
\ti{P}_3(s)&=\frac{cs((3s-1)(1+2s)(1+4s)-2cs(5+21s+10s^2-12s^3)-12c^2s^2(1+2s))}{(1+s)(1+2s)(1+3s)(1+4s)(1+6s)}.
\end{align*}
\end{Theorem}


This is a preprint of the Work accepted for publication in \href{https://www.pleiades.online/ru/journal/funan/}{Functional Analysis and Its Applications} \copyright 2020, P.Filippova 
\end{document}